\documentclass[12pt]{article}

\usepackage{amsmath}
\usepackage{amssymb}
\usepackage{amsthm}

\usepackage{ifthen}
\usepackage[dvips]{epsfig}
\usepackage{exscale}
\usepackage[dvips]{color}
\usepackage{latexsym}

\usepackage{epsf}
\usepackage{amsfonts,amssymb}
\usepackage{amsthm}

\usepackage{epsf}
\usepackage{ifthen}
\usepackage[dvips]{epsfig}
\usepackage{exscale}
\usepackage[dvips]{color}
\usepackage{latexsym}
\usepackage{amsmath}
\usepackage{amssymb}
\usepackage{makeidx}

\newtheorem{Remark}{Remark}[section]
\newtheorem{teo}{Theorem}[section]
\newtheorem{cor}{Corollary}[section]
\newtheorem{prop}{Proposition}[section]
\newtheorem{assumption}{Assumption}[section]
\newtheorem{lemma}{Lemma}[section]
\newcommand{\be}{\begin{equation}}
\newcommand{\ee}{\end{equation}}
\newcommand{\R}{\mathbb{R}}

\begin{document}
\title{A Law of Large Numbers for an Interacting Particle System with Confining Potential}
\author{Matteo Ortisi\footnote{Dept. of Mathematics, Univ. of Milano, matteoortisi@hotmail.it}}


\maketitle 

\abstract{In this paper we consider an interacting particle system modeled as a system of $N$
stochastic differential equations driven by Brownian motions with a
drift term including a confining potential acting on each particle,
and an interaction potential modeling the interaction among all the
particles of the system. The limiting behavior as the size $N$ grows
to infinity is achieved as a law of large numbers for the empirical
process associated with the interacting particle system.}

\section*{Introduction}

\medskip

We consider a system of $N(\in\mathbb{N}-\{0\})$ particles. From the Lagrangian point of view, the system is
described by $N$ random variable, $X_N^k (t) \in \mathbb R^d, t
\ge 0, k=1,\ldots, N$, so that  $\{X_N^k (t), t\in {\mathbb R}_+
\}$ is a stochastic process in the state space $({\mathbb R}^d,
\mathcal{B}_{ \mathbb R^d})$, on a common probability space
$(\Omega, \mathcal F, P)$. $X_N^k$ may describe the state of the
$k$-th particle, e.g. its position. We consider the case of
contiuous time evolution, i.e. the time evolution is described by
a system of stochastic differential equations (EDSs) with additive
noise

\begin{eqnarray}
dX_N^k (t) &=& \left[f(X_N^k(t))+ F_N \left[ X_N(t) \right] \left(
X_N^k (t) \right) \right]dt \nonumber\\
&+& \sigma_N dW^k(t), \quad k=1,\ldots, N.
\label{general_eds}
\end{eqnarray}

In equation \eqref{general_eds}  the process $\{ W^{k},
k=1,\ldots\}$ is a family of independent standard Wiener
processes, $f: \mathbb R^d+ \rightarrow \mathbb R$, and  the
functional $F_N$ is defined on ${\cal{M_P}}(\mathbb R^d)$, the
space of all probability measures on $\mathbb R^d$, and depends on
the empirical measure
\begin{equation}
X_{N}(t)=\frac{1}{N} \sum_{k=1}^{N} \epsilon_{X_{N}^{k}(t)}\in
\mathcal M_P(\mathbb R^d). \label{empirical_measure}
\end{equation}

By the empirical measure \eqref{empirical_measure} we
describe the system by an Eulerian  approach: the collective
behavior of the discrete (in the number of particles) system, may
be given in terms of the spatial distribution of particles at time
$t$.

Correspondingly, the measure valued process
\be
\label{empirical_process}
X_N:t\in\mathbb{R}_+\rightarrow X_N(t)=\frac{1}{N}\sum_{k=1}^N\epsilon_{X_N^k(t)}
\ee
is called the {\it empirical process} of the system for a population size $N$.
The trajectories are random elements of $C([0,T],\mathcal{M_P}(\mathbb{R}^d))$, so that the distributions $\mathcal{L}(X_N)$ of those processes can be considered as elements of $\mathcal{M_P}(C([0,T],\mathcal{M_P}(\mathbb{R}^d)))$.

Equation \eqref{general_eds}  might describe a system of
$N$ individuals whose movement is due to a stochastic individual
component coupled with an interaction term, and  an (individual)
advection term.

\paragraph{The individual dynamics}
The advection term $f: \mathbb R^d+ \rightarrow \mathbb R$ may
describe the individual dynamics of a particle, which may depend
on external information. Indeed we consider the following form for
$f$
\begin{eqnarray}
f (x)=-\gamma_1\nabla U(x), \label{potenzial_U}
\end{eqnarray}

where $\gamma_1\in \mathbb R_+$, and the potential $U:\mathbb
R^d\rightarrow \mathbb R_+ \in C^2(\mathbb R^d)$ is a non negative
smooth even function. From the modelling point of you the
transport term \eqref{potenzial_U} mean to be ``confining"
potential: there are some external information coming from the
environment which attracts the particle along the flow due to $U$.

\paragraph{The interaction dynamics}
$F_N$ it depends  on the relative location of the specific
individual $X_N^k(t)$ with respect to the other individuals, via
on the empirical measure of the whole system. The interaction we
consider is due to different phenomena: aggregation and repulsion.
These two different forces compete  but act at different scales.

Aggregation act at macroscale and is modelled by  a McKean-Vlasov
interaction kernel
$$
G_a : \mathbb R^d \longrightarrow \mathbb R_+.
$$

The interaction of the particle located in $X_N^k(t)$ at time $t$
with the others is described by a ``generalized'' gradient
operator as discussed in \cite{morale:MAVKDM,morale:JMB1} acting
on the empirical measure
\begin{equation}
\left(\nabla G_a * X_N(t)\right)\left(X_N^k(t)\right).
\label{F1}
\end{equation}

Repulsion acts at  mesoscale;  the mesoscale is introduced as in
\cite{morale:Int_APPL_MATH,morale:JMB1,oel2} via a rescaling of a
referring kernel $V_1$

\begin{equation}
V_N (z)= N^\beta V_1 (N^{\beta / d} z), \quad \beta \in (0,1)
\label{kernel_rescaling}
\end{equation}

The repelling force exerted on the $k$-th (out of $N$) single
particle located at $X_N^k(t)$ dis given by
\begin{eqnarray}
-\sum_{i=1}^{N}  N^{\beta - 1} \nabla\,V_1 \left(N^ {\beta / d}
\left(X_N^k (t) - X_N^i (t)\right)\right) &=& -(\nabla V_N *X_N
(t))(X_N^k(t)) \label{F2}
\end{eqnarray}

From \eqref{F2} it is clear how the choice of $\beta$
determines the range and the strength of the influence of
neighboring particles; indeed, any particle interacts (repelling)
with $O\left(N^{1-\beta}\right)$ other particles in a small volume
$O\left(N ^{-\beta}\right)$.

From \eqref{F1} and \eqref{F2}, the advection
interaction term $F[X_N]$ is given by
\begin{eqnarray}
F[X_N](x) &=& \gamma_2 \left(\nabla (G_a- V_N) *X_N (t)\right)(x),
\label{interaction_term}
\end{eqnarray}
with $\gamma_2\in \mathbb R_+$.

\paragraph{The stochasticity}

The stochastic component in equation \eqref{general_eds}
may describe both the lack of information we have about the
environment or the particle itself and the need of each particle
to interact with the others, so that they move randomly with a
mean free path $\sigma_N$ (depending on $N$) unless they meet other particles and
interact.

\medskip

By \eqref{general_eds}, \eqref{potenzial_U},
and\eqref{interaction_term} the system we study is the
following
\begin{eqnarray}
\label{system_complete}
dX_N^k(t)&=&-\left[\gamma_1\nabla U(X_N^k(t))+
\gamma_2\nabla\left(\left( V_N - G_a\right)*X_N \right)(X_N^k(t))\right]dt\nonumber\\
&& +\sigma_NdW^k(t), \quad \quad k=1,\ldots,N.
\end{eqnarray}
In the case $\gamma_1=0$, the advection is due only to the
interaction and the system become the following

\begin{eqnarray}
\label{system_VK_DM}
dX_N^k(t)&=&(\nabla G_a * X_N(t))(X_N^k(t))-(\nabla V_N * X_N(t))(X_N^k(t))dt\nonumber\\
&& +\sigma_NdW^k(t),\quad \quad k=1,\ldots,N.
\end{eqnarray}

\medskip

In previous papers
\cite{morale:MAVKDM,morale:Int_APPL_MATH,morale:JMB1} the authors
has focused their attention on the time evolution of the system
\eqref{system_VK_DM}. In particular they have analyzed the
convergence of the system as the number of particles $N$ increases
to infinity. In \cite{morale:Int_APPL_MATH,morale:JMB1} a ''law of
large numbers'' is presented while in \cite{morale:MAVKDM} the authors
have studied the existence and uniqueness of the solution of the
PDE describing the time evolution of the limit system.

In this work we focus our attention on the system
\eqref{system_complete} and extend to this case the results obtained in \cite{morale:Int_APPL_MATH,morale:JMB1}.

\section*{Notations and Hypotheses}

For some topological space $S$ we denote by $C_b^m(S,\R^d)$ the space of $m$-times differentiable $\R^d$-valued functions on $S$ with continuous bounded derivatives; $C_b^m(S,\R)$ is abbreviated with $C_b^m(S)$ and $C_b^0(S,\R^d)$ with $C_b(S,\R^d)$. $C_b^m(S,\R^d)$ is equipped with the supremum norm.
On $\R^d\times\R^d$, $( \cdot,\cdot )$ denotes the usual scalar
product.

$\mathcal{M_P}(S)$ is the space of probability measures on $S$. This space is equipped with the usual weak topology.
On the space $\mathcal{M_P}(\mathbb{R}^d)$ the weak topology is generated by the complete metric
$$
\vert\vert \mu-\nu\vert\vert_1=\sup_{f\in\mathcal{H}_1}(\langle \mu,f\rangle-\langle \nu,f\rangle),
$$
where
$$
\langle \mu,f\rangle=\int_{\mathbb{R}^d}f(x)\mu(dx)\quad f\in C_b(\mathbb{R}^d),
$$
and
$$
\mathcal{H}_1=\left\{ f\in C_b(\mathbb{R}^d)~:~\sup_{x\in\mathbb{R}^d}\vert f(x)\vert\leq 1,~\sup_{x,y\in\mathbb{R}^d,x\neq y}\frac{\vert f(x)-f(y)\vert}{\vert x-y\vert}\leq 1
\right\}.
$$
The metric $d_1(\mu,\nu):=\vert\vert \mu-\nu\vert\vert_1$ is also well known as bounded Lipschitz metric.

For any $S$-valued random variable $Y$ we denote by $\mathcal{L}(Y)\in\mathcal{M_P}(S)$ its distribution.

For some $T\in(0,\infty)$, $C([0,T],$ $\mathcal{M_P}(\mathbb{R}^d))$ is the space of all continuous functions $f=f(t),~0\leq t\leq T$ from $[0,T]$ to $\mathcal{M_P}(\mathbb{R}^d)$, equipped with the metric
$$
\rho(f,g)=\sup_{0\leq t\leq T}\vert\vert f(t)-g(t)\vert\vert_1.
$$
For $f\in L^2(\mathbb{R}^d)$ we denote by
$$
\tilde{f}(\lambda)=\lim_{a\to\infty}\left(\frac{1}{2\pi}\right)^{d/2}\int_{\{\vert x\vert\leq a\}}e^{i\lambda x}f(x)dx
$$
its Fourier transform.

In connection with Fourier transforms we shall use the relations
\begin{equation}
\label{eq_2}
\int_{\mathbb{R}^d}f(x)\overline{g(x)}dx=\int_{\mathbb{R}^d}\tilde{f}(\lambda)\overline{\tilde{g}(\lambda)}d\lambda\quad f,g\in L^2(\mathbb{R}^d),
\end{equation}
\begin{equation}
\label{eq_3}
\widetilde{f*g}(\lambda)=(2\pi)^{d/2}\tilde{f}(\lambda)\tilde{g}(\lambda)\quad f,g\in L^2(\mathbb{R}^d),
\end{equation}
\begin{equation}
\label{eq_4}
\widetilde{\nabla f}(\lambda)=-i\lambda\tilde{f}(\lambda)\quad f\in W_2^1(\mathbb{R}^d);
\end{equation}
where
$$
W_2^1(\mathbb{R}^d)=\{f\in L^2(\mathbb{R}^d)~:~\int_{\mathbb{R}^d}(1+\vert \lambda\vert^2)\vert \tilde{f}(\lambda)\vert^2d\lambda=\vert\vert f\vert\vert_2^2+\vert\vert \nabla f\vert\vert_2^2<\infty\}.
$$
Positive constants throughout the thesis are denoted by $c_1,c_2,\ldots$; if a constant depends on a quantity k, we denote it with $c(k)$.

By defining
\begin{equation}
\label{eq_1.8}
g_N(x,t)=(V_N*X_N(t))(x),
\end{equation}
\begin{equation}
\nabla g_N(x,t)=(\nabla V_N*X_N(t))(x),
\end{equation}
\begin{equation}
h_N(x,t)=(W_N*X_N(t))(x),
\end{equation}
where
\begin{equation}
\label{cond_3.4}
W_N(x)=\chi_N^dW_1(\chi_Nx)
\end{equation}
and $W_1$ is a symmetric probability density defined on $\mathbb{R}^d$, that is
\be
\label{eq_2.5_bis}
W_1(x)=W_1(-x),
\ee
equation (\ref{system_complete}) becomes
\begin{eqnarray*}
dX_N^k(t)&=&-\gamma_1\nabla U(X_N^k(t))+\gamma_2\left[(\nabla G_a*X_N(t))(X_N^k(t))-\nabla g_N(X_N^k(t),t)\right]dt\\
&+&\sigma_NdW^k(t), \quad k=1,\ldots,N.
\end{eqnarray*}
We consider the further assumptions:

\begin{equation}
\label{eq_10}
V_1(x)=(W_1*W_1)(x)=\int_{\mathbb{R}^d}W_1(x-y)W_1(y)dy,
\end{equation}

and
\begin{equation}
\label{eq_17}
W_1\in W_2^1(\R^d) {\rm \quad i.e. \quad } \int_{\mathbb{R}^d}(1+\vert \lambda\vert^2)\vert \widetilde{W_1}(\lambda)\vert^2d\lambda<\infty.
\end{equation}

About the initial condition we suppose that
\begin{equation}
\label{eq_16}
\sup_{N\in\mathbb{N}}\mathbb{E}\left[\int_{\mathbb{R}^d}\vert x\vert X_N(0)(dx)\right]<\infty,
\end{equation}

\begin{equation}
\label{eq_18}
\sup_{N\in\mathbb{N}}\mathbb{E}\left[\int_{\mathbb{R}^d}\vert h_N(x,0)\vert^2dx\right]=\sup_{N\in\mathbb{N}}\mathbb{E}\left[\vert\vert h_N(\cdot,0)\vert\vert_2^2\right]<\infty,
\end{equation}

\be
\label{cond_19}
\forall ~N\geq 1,~~\mathbf{X}_N(0)=(X^1(0),\ldots,X^N(0)) \mathrm{~is~independent~of~} W^k, k=1,\ldots,N
\ee

The assumptions about the interaction potential are the following:
\begin{itemize}
\item[]\begin{equation}
\label{eq_1.7}
G_a \in C_b^2(\mathbb{R}^d,\mathbb{R}_+)
\end{equation} is a symmetric function on $\R^d$ supposed to be independent on $N$,
\item[]$V_N$ is supposed to be of the form $V_N(x)=\chi_N^dV_1(\chi_Nx)$, where
\begin{equation}
V_1\in C_b^2(\mathbb{R}^d,\mathbb{R}_+) \mathrm{~is ~a ~symmetric ~probability ~density ~on~} \mathbb{R}^d,
\end{equation}
\item[]\begin{equation}
\chi_N=N^{\beta/d}, \quad \beta\in(0,1).
\end{equation}
\end{itemize}
It is clear that
$$
\lim_{N\to+\infty}V_N=\delta_0,
$$
where $\delta_0$ is the Dirac delta function.

\be \label{cond_20} \forall ~N\geq
1,~~\mathbb{E}\{\vert\mathbf{X}_N(0)\vert^2\}<+\infty. \ee Let the
confining potential be such that \be \label{cond_21} U\in
C_b^1(\R^d,\R_+)\cap C^2(\R^d,\R_+). \ee Then we consider the
following possible assumptions on the parameter $\beta$:
\begin{itemize}
\item[a)]\begin{eqnarray}
\label{sigma_1}
&\beta&\in\left(0,\frac{d}{d+2}\right), \lim_{N\to\infty}\sigma_N=0 \quad {\rm or}\nonumber\\  &\beta&\in\left[\frac{d}{d+2},1\right), \lim_{N\to +\infty}\sigma_NN^{\beta(d+2)/d-1}=0,
\end{eqnarray}
\item[b)]\be
\label{sigma_2}
\beta\in\left(0,\frac{d}{d+2}\right), \lim_{N\to\infty}\sigma_N=\sigma_{\infty}>0.
\ee
\end{itemize}

\section*{A Law of Large Numbers}
In this section we derive a law of large
numbers for the measure valued process $X_N$ defined by
(\ref{empirical_measure}) and (\ref{empirical_process}), in the
case of boundedness properties of the confining gradient $U$. In
particular by following the approaches proposed in \cite{mor} and \cite{oel}, we prove the existence of the limit measure 
for the sequence $\{\mathcal{L}(X_N)\}_{N\in\mathbb{N}}$ of
distributions of $\{X_N\}_{N\in\mathbb{N}}$.

We consider both the unviscous case, that is
$\lim_{N\to\infty}\sigma_N=0$, and the viscous case, i.e.
$\lim_{N\to\infty}\sigma_N=\sigma_{\infty}>0$.

The procedure may be divided into three steps:
\begin{itemize}
\item[i)]relative compactness of the sequence $\mathcal{L}(X_N), N\in\mathbb{N}$, which corresponds to an existence result of the limit $\mathcal{L}(X)$;
\item[ii)]regularity of the possible limits: we show that the possible limits $\{X(t), t\in[0,T]\}$ are absolutely continuous with respect to the Lebesgue measure for almost all $t\in [0,T]$ $\mathbb{P}-a.s.$;
\item[iii)]identification of the dynamics of the limit process: all possible limits are shown to be solution of a certain deterministic equation that we assume to have a unique solution.
\end{itemize}
In the case $\lim_{N\to\infty}\sigma_N=0$ we guess the limit
dynamics and show that it is the weak limit of $\{X_N(t), t\in
[0,T]\}$.

\subsection*{Relative Compactness}
The first step toward proving a law of large number for a measure-valued process is
to obtain a relative compactness result for the sequence of empirical measure's distribution laws $\{\mathcal{L}(X_N)\}_{N\in\mathbb{N}}$ associated to the system of stochastic differential equations.

\begin{teo}
\label{teo_1.2}
If either (\ref{sigma_1}) or (\ref{sigma_2}) holds, under conditions (\ref{eq_1.8})-(\ref{cond_21}), the sequence $\{\mathcal{L}(X_N)\}_{N\in\mathbb{N}}$ of distributions of the processes $\{X_N(t),0\leq t\leq T\}$ associated to the system of stochastic differential equations (\ref{system_complete}) is relatively compact in the space $\mathcal{M_P}(C([0,T],\mathcal{M_P}(\mathbb{R}^d)))$.
\end{teo}

We consider first some preliminary results regarding the martingale properties of some processes. Up to now we suppose
that all the hypotheses of Theorem \ref{teo_1.2} are satisfied. We
remark that all the results are valid also in the case $\gamma_1=0$.

For the seek of simplifying the notations, in the following calculations we set $\gamma_1=\gamma_2=1$.

Let
\begin{eqnarray}
\label{A_N}
A_N(t)&=&\int_0^t\langle X_N(s),2(\vert \nabla g_N(\cdot,u)\vert^2-\nabla g_N(\cdot,u)(-\nabla U(\cdot)+(\nabla G_a*X_N(u))(\cdot))\nonumber\\
&&+\vert -\nabla U(\cdot)+(\nabla G_a*X_N(u))(\cdot)\vert^2)\rangle
+\sigma_N^2\vert\vert \nabla h_N(\cdot,u)\vert\vert_2^2du.
\end{eqnarray}

\begin{lemma}
\label{lemma_4.1}
The process
\begin{eqnarray*}
M_N(t)&=&\vert\vert h_N(\cdot,t)\vert\vert_2^2+A_N(t)-\int_0^t\langle X_N(u),2\vert-\nabla U(\cdot)+(\nabla G_a*X_N(u))(\cdot)\vert^2\rangle du\\
&-&c_1\sigma_N^2 t N^{\beta (d+2)/d-1}
\end{eqnarray*}
is a martingale.
\end{lemma}
\begin{flushleft}
{\it    Proof.}
\end{flushleft}
Because of (\ref{eq_1.7}), (\ref{eq_2.5_bis}), (\ref{eq_17}) and (\ref{cond_21}),
by applying Ito's formula to
$$
\vert\vert h_N(\cdot,t)\vert\vert_2^2=\langle X_N(t),g_N(\cdot,t)\rangle=\frac{1}{N^2}
\sum_{k,l=1}^NV_N(X_N^k(t)-X_N^l(t)),
$$
one obtains
\begin{eqnarray}
\label{eq_2.15}
\mathbb{E}\left[\vert\vert h_N(\cdot,t)\vert\vert_2^2\vert\mathcal{F}_s\right]
&=&\vert\vert h_N(\cdot,s)\vert\vert_2^2\nonumber\\
&+&\mathbb{E}\left[\frac{2}{N}\int_s^t\sum_{m=1}^N\langle X_N(u),\nabla G_a(\cdot-X_N^m(u))\cdot\nabla g_N(\cdot,u)\rangle du\right.\nonumber\\
&+&\frac{2}{N}\int_s^t\sum_{m=1}^N\langle X_N(u),-\nabla U(\cdot)\cdot\nabla g_N(\cdot,u)\rangle du\nonumber\\
&-&\frac{2}{N}\int_s^t\sum_{m=1}^N\langle X_N(u),(\nabla V_N(\cdot-X_N^m(u)))\cdot\nabla g_N(\cdot,u)\rangle du\nonumber\\
&+&\left.\frac{\sigma_N^2}{N^2}\int_s^t\sum_{k,m=1, ~k\neq m}^N\Delta V_N(X_N^k(u)-X_N^m(u))du\vert\mathcal{F}_s\right]\nonumber\\
&=&\vert\vert h_N(\cdot,s)\vert\vert_2^2\nonumber\\
&-&\mathbb{E}\left[2\int_s^t\langle X_N(u),\vert \nabla g_N(\cdot,u)\vert^2\rangle du \right.\nonumber\\
&-&2\int_s^t\langle X_N(u),\nabla g_N(\cdot,u)\cdot(-\nabla U(\cdot)+(\nabla G_a*X_N(u))(\cdot))\rangle du\nonumber\\
&+&\left. \sigma_N^2\int_s^t\vert\vert \nabla h_N(\cdot,u)\vert\vert_2^2 du\vert\mathcal{F}_s\right] +\frac{\sigma_N^2 (t-s)}{N}\Delta V_N(0).
\end{eqnarray}

By assumption (\ref{eq_17}),
\begin{eqnarray}
\Delta V_1(0)&=&\Delta (W_1*W_1)(0)=(\nabla W_1*\nabla W_1)(0)\nonumber\\
&=&\int_{\mathbb{R}^d}\nabla W_1(-y)\nabla W_1(y)dy\nonumber\\
&\leq&\left(\int_{\mathbb{R}^d}\vert \nabla W_1(-y)\vert^2dy\right)^{1/2}\left(\int_{\mathbb{R}^d}\vert \nabla W_1(y)\vert^2dy\right)^{1/2}\nonumber\\
&<&\infty.
\end{eqnarray}

As a consequence
\begin{eqnarray}
\vert\Delta V_N(x)\vert=N^{\beta}N^{2\beta/d}\vert\Delta V_1(N^{\beta/d} x)\vert.
\end{eqnarray}

So
\begin{equation}
\frac{\sigma_N^2 (t-s)}{N}\Delta V_N(0)=c_1 \sigma_N^2(t-s) N^{\beta(d+2)/d-1},
\end{equation}
and the thesis follows.

\hspace{14 cm}$\Box$

\begin{Remark}
$A_N(t)$ is not negative; indeed in general
$$
a^2-ab+b^2\geq \frac{a^2}{2}-ab+\frac{b^2}{2}=(\frac{a}{\sqrt{2}}-\frac{b}{\sqrt{2}})^2\geq0.
$$
\end{Remark}

Let now consider a special class of test functions, i.e. positive function $\phi \in C^2(\mathbb{R}^d)$ such that
\be
\label{cond_phi}
\phi(x)=\vert x\vert\quad \mathrm{for~}\vert x\vert\geq 1\quad \mathrm{and}\quad \vert\vert \nabla\phi\vert\vert_{\infty}+\vert\vert \Delta\phi\vert\vert_{\infty}<\infty.
\ee
\begin{lemma}
\label{lemma_4.2}
For $0\leq s<t\leq T$ and $\phi\in C^2(\R^d)$ such that (\ref{cond_phi}) holds
\begin{equation}
\label{eq_1.38}
\langle X_N(t),\phi\rangle+c_2A_N(t)+c_3t \quad \mathit{is ~a~ submartingale}
\end{equation}
and
\begin{equation}
\label{eq_36}
\langle X_N(t),\phi\rangle-c_2A_N(t)-c_3t \quad \mathit{is ~a~ supermartingale},
\end{equation}
with $c_2, c_3\in\R_+$.
\end{lemma}
\begin{flushleft}
{\it    Proof.}
\end{flushleft}
By applying Ito's Formula to $\langle X_N(t),\phi\rangle$,
\begin{eqnarray}
\label{eq_1.40}
&&\mathbb{E}\left[\langle X_N(t),\phi\rangle\vert\mathcal{F}_s\right]\nonumber\\
&&=\langle X_N(s),\phi\rangle\nonumber\\
&&+\mathbb{E}\left[\int_s^t\langle X_N(u),(-\nabla U(\cdot)+(\nabla G_a*X_N(u))(\cdot)-\nabla g_N(\cdot,u))\nabla\phi\right.\nonumber\\
&&\left.+\frac{\sigma_N^2}{2}\Delta\phi\rangle du\vert\mathcal{F}_s\right]\nonumber\\
&&\geq\langle X_N(s),\phi\rangle\nonumber\\
&&-c_2\mathbb{E}\left[\int_s^t\langle X_N(u),-\nabla U(\cdot)+(\nabla G_a*X_N(u))(\cdot)-\nabla g_N(\cdot,u)+1\rangle du\vert\mathcal{F}_s\right].\nonumber\\
\end{eqnarray}
Since
\begin{eqnarray*}
0&\leq&\langle X_N(u),\vert -\nabla U(\cdot)+(\nabla G_a*X_N(u))(\cdot)-\nabla g_N(\cdot,u)-1\vert^2\rangle\\
&=&\langle X_N(u),\vert -\nabla U(\cdot)+(\nabla G_a*X_N(u))(\cdot)-\nabla g_N(\cdot,u)\vert^2\\
&-&2(-\nabla U(\cdot)+(\nabla G_a*X_N(u))(\cdot)-\nabla g_N(\cdot,u))+1\rangle,
\end{eqnarray*}
\begin{eqnarray*}
&&2\langle X_N(u),-\nabla U(\cdot)+(\nabla G_a*X_N(u))(\cdot)-\nabla g_N(\cdot,u)\rangle\\
&&\leq \langle X_N(u),\vert -\nabla U(\cdot)+(\nabla G_a*X_N(u))(\cdot)-\nabla g_N(\cdot,u)\vert^2+1\rangle,
\end{eqnarray*}
and therefore
\begin{eqnarray*}
&&\langle X_N(u),-\nabla U(\cdot)+(\nabla G_a*X_N(u))(\cdot)-\nabla g_N(\cdot,u)\rangle\\
&&\leq \langle X_N(u),\vert -\nabla U(\cdot)+(\nabla G_a*X_N(u))(\cdot)-\nabla g_N(\cdot,u)\vert^2+1\rangle.
\end{eqnarray*}
This implies that (\ref{eq_1.40}) is greater than or equal to
\begin{eqnarray}
\label{eq_38}
&&\langle X_N(s),\phi\rangle
-c_2\mathbb{E}\left[\int_s^t\langle X_N(u),\vert -\nabla U(\cdot)+(\nabla G_a*X_N(u))(\cdot)-\nabla g_N(\cdot,u)\vert^2\right.\nonumber\\
&&\left.+2\rangle du\vert\mathcal{F}_s\right]\nonumber\\
&&\geq \langle X_N(s),\phi\rangle
-c_3\mathbb{E}\left[\int_s^t\langle X_N(u),\vert -\nabla U(\cdot)+(\nabla G_a*X_N(u))(\cdot)-\nabla g_N(\cdot,u)\vert^2\right.\nonumber\\
&&\left.+1\rangle du\vert\mathcal{F}_s\right]\nonumber\\
&&\geq \langle X_N(s),\phi\rangle
-c_3\mathbb{E}\left[\int_s^t\langle X_N(u),2(\vert -\nabla U(\cdot)+(\nabla G_a*X_N(u))(\cdot)\vert^2\right.\nonumber\\
&&\left.-(-\nabla U(\cdot)+(\nabla G_a*X_N(u))(\cdot))\nabla g_N(\cdot,u)
+\vert\nabla g_N(\cdot,u)\vert^2)\rangle+\sigma_N^2\vert\vert \nabla h_N(\cdot,u)\vert\vert_2^2 du\right.\nonumber\\
&&\left.+\int_s^tdu\vert\mathcal{F}_s\right]\nonumber\\
&&=\langle X_N(s),\phi\rangle-c_3\mathbb{E}\left[A_N(t)-A_N(s)+t-s\vert\mathcal{F}_s\right]\nonumber\\
&&=\langle X_N(s),\phi\rangle-c_3\mathbb{E}\left[A_N(t)+t\vert\mathcal{F}_s\right]+c_3A_N(s)+c_3s.
\end{eqnarray}
Hence,
$$
\mathbb{E}\left[\langle X_N(t),\phi\rangle+c_3A_N(t)+c_3t\vert\mathcal{F}_s\right]\geq\langle X_N(s),\phi\rangle+c_3A_N(s)+c_3s
$$
and (\ref{eq_1.38}) follows.

In a completely analogous way (with $+c_3$ instead of $-c_3$), we obtain the property (\ref{eq_36}).

\hspace{14 cm}$\Box$

Let us define the sequence of stopped processes
$X_{N,k}(t)=X_N(t\wedge\tau_N^k),~0\leq t\leq T,~N\in\mathbb{N}$, and $k>0$ fixed, where $\tau_{N}^k$ is defined by
\begin{eqnarray}
\label{tau}
\!\!\!\!\!\!\!\!\!\!\!\!\tau_N^k&=&\inf\{t\geq 0~:~\vert\vert h_N(\cdot,t)\vert\vert_2^2+A_N(t)\nonumber\\
&&-\int_0^t\langle X_N(u),2\vert-\nabla U(\cdot)+(\nabla G_a*X_N(u))(\cdot)\vert^2\rangle du >k\},
\end{eqnarray}
where $k\in\R+$.

We consider a slight modification of the more general characterization of the relative compactness
 by Ethier-Kurtz (see \cite{kurtz}, theorem 8.2). In particular we prove confining compactness property for the process $\{X_{N,k}(t)\}$ and then the boundedness of small variations of the process.
\begin{prop}
\label{prop_1.1}
For any $\epsilon>0$ there exists a compact $K_{\epsilon}^k$ in $(\mathcal{M_P}(\mathbb{R}^d),\vert\vert\cdot\vert\vert_1)$ such that
$$
\inf_{N\in\mathbb{N}}\mathbb{P}\{X_{N,k}(t)\in K_{\epsilon}^k,~\forall t\in[0,T]\}\geq 1-\epsilon.
$$
\end{prop}
\begin{flushleft}
{\it    Proof.}
\end{flushleft}
Let $B_{\lambda}^c=\{x\in\mathbb{R}^d~:~\vert x\vert>\lambda, ~\lambda>1\}$;
\begin{eqnarray}
\langle X_{N,k}(t),\phi\rangle&=&\int_{\mathbb{R}^d}\phi(x)X_{N,k}(t)(dx)\geq\int_{\{x:\vert x\vert>\lambda,~\lambda>1\}}\phi(x)X_{N,k}(t)(dx)\nonumber\\
&=&\int_{\{x:\vert x\vert>\lambda,~\lambda>1\}}\vert x\vert X_{N,k}(t)(dx)\geq \lambda\int_{\{x:\vert x\vert>\lambda,~\lambda>1\}} X_{N,k}(t)(dx)\nonumber\\
&=&\lambda\langle X_{N,k}(t),\mathbf{1}_{B_{\lambda}^c}\rangle;\nonumber
\end{eqnarray}
if $\langle X_{N,k}(t),\mathbf{1}_{B_{\lambda}^c}\rangle>\delta$, then $\langle X_{N,k}(t),\phi\rangle>\lambda\delta$ and therefore
\begin{eqnarray}
\label{eq_39}
&&\mathbb{P}\left\{\sup_{t\leq T} \langle X_{N,k}(t),\mathbf{1}_{B_{\lambda}^c}\rangle>\delta\right\}\leq\mathbb{P}\left\{\sup_{t\leq T}\langle X_{N,k}(t),\phi\rangle>\lambda\delta\right\}\nonumber\\
&&\leq\mathbb{P}\left\{\sup_{t\leq T}\langle X_{N,k}(t),\phi\rangle+c_3A_N(t\wedge\tau_N^k)+(t\wedge\tau_N^k)c_3>\lambda\delta\right\}.
\end{eqnarray}
By (\ref{eq_1.38}) in Lemma \ref{lemma_4.2} and by Doob's Inequality, (\ref{eq_39}) is less than or equal to
\begin{eqnarray}
\label{eq_40}
&&\frac{1}{\lambda\delta}\left(\mathbb{E}\left[\langle X_{N,k}(T),\phi\rangle\right]+c_3\mathbb{E}\left[A_N(T\wedge\tau_N^k)\right]+c_3\mathbb{E}\left[T\wedge\tau_N^k\right]\right)\nonumber\\
&&\leq\frac{1}{\lambda\delta}\left(\mathbb{E}\left[\mathbb{E}[\langle X_{N,k}(T),\phi\rangle\vert\mathcal{F}_0]\right]+c_3\mathbb{E}\left[A_N(T\wedge\tau_N^k)\right]+c_3\mathbb{E}\left[T\wedge\tau_N^k\right]\right)\nonumber\\
&&\stackrel{(\ref{eq_36})}{\leq}\frac{1}{\lambda\delta}\left(\mathbb{E}\left[\langle X_{N,k}(0),\phi\rangle+c_3\mathbb{E}[A_N(T\wedge\tau_N^k)\vert\mathcal{F}_0]+c_3\mathbb{E}[T\wedge\tau_N^k\vert\mathcal{F}_0]\right.\right.\nonumber\\
&&\left.\left.-c_3A_N(0\wedge\tau_N^k)-(0\wedge\tau_N^k)c_3+c_3A_N(T\wedge\tau_N^k)+(T\wedge\tau_N^k)c_3\right]\right).
\end{eqnarray}
By the definition of $\tau_N^k$ and since $\nabla G_a, \nabla U \in C_b(\mathbb{R}^d,\mathbb{R}_+)$,
\begin{eqnarray}
A_N(T\wedge\tau_N^k)&<&k-\vert\vert h_N(\cdot,t\wedge\tau_N^k)\vert\vert_2^2\nonumber\\
&&+\int_0^{T\wedge\tau_N^k}\langle X_N(u),2\vert -\nabla U(\cdot)+(\nabla G_a*X_N(u))(\cdot)\vert^2\rangle du\nonumber\\
&<&k+\int_0^{T\wedge\tau_N^k}\langle X_N(u),2\vert -\nabla U(\cdot)+(\nabla G_a*X_N(u))(\cdot)\vert^2\rangle du\nonumber\\
&<&k+c_4T;
\end{eqnarray}
it follows that (\ref{eq_40}) is less than or equal to
\begin{eqnarray}
\frac{1}{\lambda\delta}\left(\mathbb{E}\left[\langle X_{N,k}(0),\phi\rangle+c_3(k+c_4T)+c_3T+c_3(k+c_4T)+c_3T\right]\right)\stackrel{(\ref{eq_16})}{\leq}\frac{c_5(k,T)}{\lambda\delta}.
\end{eqnarray}
As a consequence,
\begin{equation}
\label{eq_1.46}
\mathbb{P}\left\{\sup_{t\leq T} \langle X_{N,k}(t),\mathbf{1}_{B_{\lambda}^c}\rangle>\delta\right\}\leq\frac{c_5(k,T)}{\lambda\delta}.
\end{equation}
Let us now take $\epsilon>0$ and two sequences $\mu_i$ and $\delta_i$ of positive numbers such that $\sum_{i=1}^{\infty}\mu_i=\epsilon$ and $\delta_i\searrow 0$. Let $\lambda_i=\frac{c_5(k,T)}{\mu_i\delta_i}\to\infty$. Then (\ref{eq_1.46}) yelds
\begin{eqnarray}
\label{eq_44}
\mathbb{P}\left\{\sup_{t\leq T}\langle X_{N,k}(t),\mathbf{1}_{B_{\lambda_i}^c}\rangle> \delta_i,\quad\forall i\in\mathbb{N}\right\}&\leq&\sum_{i=1}^{\infty}\mathbb{P}\left\{\sup_{t\leq T}\langle X_{N,k}(t),\mathbf{1}_{B_{\lambda_i}^c}\rangle> \delta_i\right\}\nonumber\\
&\leq&\sum_{i=1}^{\infty}\frac{c_5(k,T)}{\lambda_i\delta_i}=\sum_{i=1}^{\infty}\mu_i=\epsilon.
\end{eqnarray}
By Prohorov's Theorem, the set
$$
K_{\epsilon}^k=\{\mu\in\mathcal{M_P}(\mathbb{R}^d)~:~\langle \mu,\mathbf{1}_{B_{\lambda_i}^c}\rangle\leq \delta_i,\quad \forall i\in\mathbb{N}\}
$$
is compact in $\mathcal{M_P}(\mathbb{R}^d)$; since
$$
\mathbb{P}\left\{\sup_{t\leq T}\langle X_{N,k}(t),\mathbf{1}_{B_{\lambda_i}^c}\rangle> \delta_i,\!\forall i\in\mathbb{N}\right\}=1-\mathbb{P}\left\{\langle X_{N,k}(t),\mathbf{1}_{B_{\lambda_i}^c}\rangle\leq\delta_i,\! \forall i\in\mathbb{N},\! \forall t\in[0,T]\right\},
$$
by (\ref{eq_44}), $\forall \epsilon >0$ there exists a compact set $K_{\epsilon}^k\subset\mathcal{M_P}(\mathbb{R}^d)$ such that
$$
\inf_{N\in\mathbb{N}}\mathbb{P}\{X_{N,k}(t)\in K_{\epsilon}^k,~\forall t\in[0,T]\}\geq 1-\epsilon.
$$

\hspace{14 cm}$\Box$

Next proposition states that for little time variations we have little variations of the processes $\{X_{N,k}(t)\}$.

\begin{prop}
\label{prop_4.2}
For any $0<\delta<1$, there exists a sequence $\{\gamma_n^T(\delta)\}_{n\in\mathbb{N}}$ of non negative random variables such that
\begin{equation}
\mathbb{E}\left[\vert\vert X_{N,k}(t+\delta)-X_{N,k}(t)\vert\vert_1^4\right]\leq\mathbb{E}\left[\gamma_n^T(\delta)\right]\quad 0\leq t\leq T
\end{equation}
and
\begin{equation}
\lim_{\delta\to 0}\limsup_{n\to\infty}\mathbb{E}[\gamma_n^T(\delta)]=0.
\end{equation}
\end{prop}
\begin{flushleft}
{\it Proof.}
\end{flushleft}
\begin{eqnarray}
\label{eq_47}
&&\vert\vert X_{N,k}(t)-X_{N,k}(s)\vert\vert_1=\sup_{f\in\mathcal{H}_1}\frac{1}{N}\sum_{i=1}^N\left( f(X_{N,k}^i(t))-f(X_{N,k}^i(s))\right)\nonumber\\
&&\leq \frac{1}{N}\sum_{i=1}^N \vert X_{N,k}^i(t)-X_{N,k}^i(s)\vert\nonumber\\
&&= \frac{1}{N}\sum_{i=1}^N\vert X_{N}^i(t\wedge \tau_N^k)-X_{N}^i(s\wedge \tau_N^k)\vert\nonumber\\
&&= \frac{1}{N}\sum_{i=1}^N\left\vert \int_{s\wedge \tau_N^k}^{t\wedge \tau_N^k}-\nabla U(X_N^i(u))+(\nabla G_a*X_N(u))(X_N^i(u))-\nabla g_N(X_N^i(u),u)du\right.\nonumber\\
&&\left.+\sigma_N(W^i(t\wedge \tau_N^k)-W^i(s\wedge \tau_N^k))\right\vert\nonumber\\
&&\leq \frac{c_7}{N}\sum_{i=1}^N\left\vert\int_{s\wedge \tau_N^k}^{t\wedge \tau_N^k}du \right\vert+\nonumber\\ &&\frac{1}{N}\sum_{i=1}^N \int_{s\wedge \tau_N^k}^{t\wedge \tau_N^k}\left\vert(\nabla G_a*X_N(u))(X_N^i(u))-\nabla g_N(X_N^i(u),u)\right\vert du\nonumber\\
&&+\frac{\sigma_N}{N}\sum_{i=1}^N\left\vert W^i(t\wedge \tau_N^k)-W^i(s\wedge \tau_N^k)\right\vert\nonumber\\
\end{eqnarray}
\begin{eqnarray}
&&=\frac{c_7}{N}\sum_{i=1}^N\left\vert\int_{s\wedge \tau_N^k}^{t\wedge \tau_N^k}du \right\vert+\int_{s\wedge \tau_N^k}^{t\wedge \tau_N^k}\langle X_N(u),\vert (\nabla G_a*X_N(u))(\cdot)-\nabla g_N(\cdot,u)\vert\rangle du\nonumber\\
&&+\frac{\sigma_N}{N}\sum_{i=1}^N\left\vert W^i(t\wedge \tau_N^k)-W^i(s\wedge \tau_N^k)\right\vert.
\end{eqnarray}
By the Cauchy-Schwartz and Jensen inequalities,
\begin{eqnarray}
\label{eq_48}
&&\int_{s\wedge \tau_N^k}^{t\wedge \tau_N^k}\langle X_N(u),\vert (\nabla G_a*X_N(u))(\cdot)-\nabla g_N(\cdot,u)\vert\rangle du\nonumber\\
&&\leq\left( \int_{s\wedge \tau_N^k}^{t\wedge \tau_N^k} du\right)^{1/2}\left(\int_{s\wedge \tau_N^k}^{t\wedge \tau_N^k}\langle X_N(u),\vert (\nabla G_a*X_N(u))(\cdot)-\nabla g_N(\cdot,u)\vert\rangle^2 du\right)^{1/2}\nonumber\\
&&\leq\left( \int_{s\wedge \tau_N^k}^{t\wedge \tau_N^k} du\right)^{1/2}\left(\int_{s\wedge \tau_N^k}^{t\wedge \tau_N^k}\langle X_N(u),\vert (\nabla G_a*X_N(u))(\cdot)-\nabla g_N(\cdot,u)\vert^2\rangle du\right)^{1/2};\nonumber\\
\end{eqnarray}
moreover, if $s\leq \tau_N^k$,
\begin{equation}
\label{eq_49}
\int_{s\wedge \tau_N^k}^{t\wedge \tau_N^k}du=\int_{s}^{t\wedge \tau_N^k}du\leq\int_{s}^{t}du=t-s
\end{equation}
and if $s> \tau_N^k$
\begin{equation}
\label{eq_50}
\int_{s\wedge \tau_N^k}^{t\wedge \tau_N^k}du=\int_{\tau_N^k}^{\tau_N^k}du=0\leq t-s.
\end{equation}
Therefore, by (\ref{eq_48}),(\ref{eq_49}) and (\ref{eq_50}), (\ref{eq_47}) is less than or equal to
\begin{eqnarray}
&&c_7\vert t-s\vert+(t-s)^{1/2}\left(\int_{0}^{t\wedge \tau_N^k}\langle X_N(u),\vert (\nabla G_a*X_N(u))(\cdot)-\nabla g_N(\cdot,u)\vert^2\rangle du\right)^{1/2}\nonumber\\
&&+\frac{\sigma_N}{N}\sum_{i=1}^N\left\vert W^i(t\wedge \tau_N^k)-W^i(s\wedge \tau_N^k)\right\vert\nonumber\\
&&\leq c_7\vert t-s\vert+(t-s)^{1/2}\left(\int_0^{t\wedge \tau_N^k}\langle X_N(u),2\vert (\nabla G_a*X_N(u))(\cdot)\vert^2\right.\nonumber\\
&&\left.-2((\nabla G_a*X_N(u))(\cdot)\nabla g_N(\cdot,u))+2\vert \nabla g_N(\cdot,u)\vert^2\rangle\right. \nonumber\\
&&\left.+\sigma_N^2\vert\vert \nabla h_N(\cdot,T\wedge\tau_N^k)\vert\vert_2^2du+\vert\vert h_N(\cdot,T\wedge \tau_N^k)\vert\vert_2^2\right)^{1/2}\nonumber\\
&&+\frac{\sigma_N}{N}\sum_{i=1}^N\left\vert W^i(t\wedge \tau_N^k)-W^i(s\wedge \tau_N^k)\right\vert\nonumber
\end{eqnarray}
\begin{eqnarray}
&&=c_7\vert t-s\vert+(t-s)^{1/2}\left(A_N(T\wedge\tau_N^k)+\vert\vert h_N(\cdot,T\wedge\tau_N^k)\vert\vert_2^2\right)^{1/2}\nonumber\\
&&+\frac{\sigma_N}{N}\sum_{i=1}^N\left\vert W^i(t\wedge \tau_N^k)-W^i(s\wedge \tau_N^k)\right\vert\nonumber\\
&&\leq c_7\vert t-s\vert+(t-s)^{1/2}\left(k+\int_0^{T\wedge\tau_N^k}\langle X_N(s),2\vert (\nabla G_a*X_N(u))(\cdot)\vert^2\rangle du\right)^{1/2}\nonumber\\
&&+\frac{\sigma_N}{N}\sum_{i=1}^N\left\vert W^i(t\wedge \tau_N^k)-W^i(s\wedge \tau_N^k)\right\vert\nonumber\\
&&\leq c_7\vert t-s\vert+(t-s)^{1/2}(k+c_6T)^{1/2}+\frac{\sigma_N}{N}\sum_{i=1}^N\left\vert W^i(t\wedge \tau_N^k)-W^i(s\wedge \tau_N^k)\right\vert.\nonumber\\
\end{eqnarray}

It follows that
\be
\vert\vert X_{N,k}(t)-X_{N,k}(s)\vert\vert_1^4\leq 2^6\left[c_7^4(t-s)^4+(t-s)^2(k+c_6T)^2+\frac{\sigma_N^4}{N^3}(t-s)^2\right].
\ee

As a consequence, for $0\leq s<t\leq T$,
$$
\mathbb{E}[\vert\vert X_{N,k}(t)-X_{N,k}(s)\vert\vert_1^4]\leq 2^6\mathbb{E}\left[c_7^4(t-s)^4+(t-s)^2(k+c_6T)^2+\frac{\sigma_N^4}{N^3}(t-s)^2\right];
$$

in particular, with $t-s=\delta$ and
$$\gamma_N^T:\delta\mapsto 2^6\left[c_7^4\delta^4+\delta^2(k+c_6T)^2+\frac{\sigma_N^4}{N^3}\delta^2\right],$$

we obtain
$$
\mathbb{E}\left[\vert\vert X_{N,k}(t)-X_{N,k}(s)\vert\vert_1^4\right]\leq\mathbb{E}[\gamma_N^T(\delta)]
$$

and
$$
\lim_{\delta\to 0}\sup_{N\in\mathbb{N}}\mathbb{E}[\gamma_N^T(\delta)]=0.
$$

\hspace{14 cm}$\Box$

\begin{prop}
\label{prop_1.3}
$\{\mathcal{L}(X_N(\cdot\wedge\tau_N^k))\}_{N\in\mathbb{N}}$, the sequence of probability laws of the processes $\{X_N((t\wedge\tau_N^k)),~0\leq t\leq T\}$ is relatively compact in $\mathcal{M_P}(C([0,T],\mathcal{M_P}(\mathbb{R}^d)))$.
\end{prop}
\begin{flushleft}
{\it    Proof.}
\end{flushleft}
It is an obvious consequence of Proposition \ref{prop_4.2} and Theorem 8.6 p.137, in \cite{kurtz}.

\hspace{14 cm}$\Box$

\begin{prop}
\label{prop_4.4}
For any $\tau$ such that $0<\tau<\infty$,
$$
\lim_{k\to\infty}\inf_{N\in\mathbb{N}}\mathbb{P}\{\tau_N^k>\tau\}=1.
$$
\end{prop}
\begin{flushleft}
{\it    Proof.}
\end{flushleft}
By Lemma \ref{lemma_4.1}, the process
$$
t\mapsto S_N(t)=\vert\vert h_N(\cdot,t)\vert\vert_2^2+A_N(t)-\int_0^t\langle X_N(u),2\vert-\nabla U(\cdot)+(\nabla G_a*X_N(u))(\cdot)\vert^2\rangle du
$$
is a submartingale.

By Doob's inequality
\begin{eqnarray}
\label{eq_45}
\mathbb{P}\left\{\sup_{t\leq \tau}S_N(t)>k\right\}\leq \frac{1}{k}\mathbb{E}[S_N(\tau)]&=&\frac{1}{k}\mathbb{E}[M_N(\tau)+\tau\sigma_N^2N^{\beta(d+2)/d-1}c_1]\nonumber\\
&=&\frac{1}{k}\mathbb{E}\left[\mathbb{E}[M_N(\tau)\vert\mathcal{F}_0]+\tau\sigma_N^2N^{\beta(d+2)/d-1}c_1\right]\nonumber\\
&=&\frac{1}{k}\mathbb{E}[M_N(0)+\tau\sigma_N^2N^{\beta(d+2)/d-1}c_1]\nonumber\\
&=&\frac{1}{k}\left(\mathbb{E}[\vert\vert h_N(\cdot,0)\vert\vert_2^2]+\tau\sigma_N^2N^{\beta(d+2)/d-1}c_1\right);\nonumber\\
\end{eqnarray}
since $\lim_{N\to\infty}\sigma_N=\sigma_{\infty}\geq 0$, by (\ref{eq_18}) and (\ref{sigma_1}) or (\ref{sigma_2}), (\ref{eq_45}) is less than or equal to $c_8(\tau)/k$, uniformly in $N$.

The thesis follows.

\hspace{14 cm}$\Box$

\begin{Remark}
\label{rem_4.2}
Proposition \ref{prop_4.4} implies that $~t\wedge\tau_N^k=t$, for any $\tau$ such that $0\leq t\leq \tau$.
\end{Remark}

\begin{flushleft}
\underline{{\it Proof of Theorem \ref{teo_1.2}.}}
\end{flushleft}
At this point Theorem \ref{teo_1.2} simply follows from Propositions \ref{prop_1.3}, \ref{prop_4.4} and Remark \ref{rem_4.2}.

\hspace{14 cm}$\Box$

Theorem \ref{teo_1.2} implies the existence of a subsequence $N_k\subset\mathbb{N}$, $N_1<N_2<\ldots$, such that the sequence $\{\mathcal{L}(X_{N_k})\}_{k\in\mathbb{N}}$ converges in $\mathcal{M_P}(C([0,T],\mathcal{M_P}(\R^d)))$ to some limit $\mathcal{L}(X)$, which is the distribution of some process $\{X(t), t\in[0,T]\}$, with trajectories in $C([0,T],\mathcal{M_P}(\R^d))$.
We discuss the uniqueness of the limit later on. By now we assume the uniqueness, so that $\{N_k\}=\mathbb{N}$.

By Skorokhod's Theorem, we are allowed to assume that $\{X_N(t), t\in [0,T]\}$ converges $\mathbb{P}$-almost surely to $\{X(t), t\in [0,T]\}$ as $N$ grows to infinity.
So, we have
\be
\label{eq_3.6}
\lim_{N\to\infty}\sup_{t\leq T}\vert\vert X_N(t)-X(t)\vert\vert_1=0\quad \mathbb{P}-a.s.
\ee

\section*{Absolute continuity of the limit}
Next proposition deals with the regularity properties of the limit measure $X(t)$. We consider the viscous case $\lim_{N\to\infty}\sigma_N=\sigma_{\infty}>0$.

\begin{prop}
\label{prop_3.2} Suppose that
$\lim_{N\to\infty}\sigma_N=\sigma_{\infty}>0$. For any $t\geq 0$,
the measure $X(t)$ is absolutely continuous with respect to
Lebesgue measure on $\R^d$ with a density $\rho\in
L^2(\R^d,\R_+)$.
\end{prop}
\begin{flushleft}
{\it Proof.}
\end{flushleft}
We begin by showing that there exists a positive function $\rho(x,t)$ such that
\be
\label{density}
\lim_{N\to\infty}\mathbb{E}\left[\int_0^T\int_{\R^d}\vert h_N(x,t)-\rho(x,t)\vert^2 dxdt\right]=0.
\ee
\begin{eqnarray}
\label{sys_3.7}
&&\lim_{N,N'\to\infty}\mathbb{E}\left[\int_0^T\int_{\R^d}\vert h_N(x,t)-h_{N'}(x,t)\vert^2 dxdt\right]\nonumber\\
&=&\lim_{N,N'\to\infty}\mathbb{E}\left[\int_0^T\int_{\R^d}\vert \widetilde{h_N}(\lambda,t)-\widetilde{h_{N'}}(\lambda,t)\vert^2 d\lambda dt\right]\nonumber\\
&\leq&\lim_{N,N'\to\infty}\mathbb{E}\left[\int_0^T\int_{\{\vert \lambda\vert\leq k\}}\vert \widetilde{h_N}(\lambda,t)-\widetilde{h_{N'}}(\lambda,t)\vert^2 d\lambda dt\right]\nonumber\\
&&+2\lim_{N,N'\to\infty}\mathbb{E}\left[\int_0^T\int_{\{\vert \lambda\vert> k\}}\vert \widetilde{h_N}(\lambda,t)\vert^2+\vert\widetilde{h_{N'}}(\lambda,t)\vert^2 d\lambda dt\right];
\end{eqnarray}
since, by (\ref{eq_17}), $\vert \widetilde{W_N}(\lambda)\vert$ is bounded and $\widetilde{h_N}(\lambda ,t)=\langle X_N(t),e^{i\lambda\cdot}\rangle\widetilde{W_N}(\lambda)$, expression (\ref{sys_3.7}) is less than or equal to
\begin{eqnarray}
\label{sys_3.8}
&&\lim_{N,N'\to\infty}\mathbb{E}\left[\int_0^T\int_{\{\vert \lambda\vert\leq k\}}\vert \langle X_N(t),e^{i\lambda\cdot}\rangle-\langle X_{N'}(t),e^{i\lambda\cdot}\rangle\vert^2 d\lambda dt\right]\nonumber\\
&&+2\lim_{N,N'\to\infty}\mathbb{E}\left[\int_0^T\int_{\{\vert \lambda\vert>k\}}\frac{\vert\lambda\vert^2}{k}\vert \widetilde{h_N}(\lambda,t)\vert^2+\frac{\vert\lambda\vert^2}{k}\vert\widetilde{h_{N'}}(\lambda,t)\vert^2 d\lambda dt\right]\nonumber\\
&\leq&\lim_{N,N'\to\infty}\mathbb{E}\left[\int_0^T\int_{\{\vert \lambda\vert\leq k\}}\vert \langle X_N(t),e^{i\lambda\cdot}\rangle-\langle X_{N'}(t),e^{i\lambda\cdot}\rangle\vert^2 d\lambda dt\right]\nonumber\\
&&+\frac{2}{k}\lim_{N,N'\to\infty}\mathbb{E}\left[\int_0^T\int_{\R^d}\vert\lambda\vert^2\vert \widetilde{h_N}(\lambda,t)\vert^2+\vert\lambda\vert^2\vert\widetilde{h_{N'}}(\lambda,t)\vert^2 d\lambda dt\right].
\end{eqnarray}
Since by (\ref{eq_3.6})
$$
\lim_{N,N'\to\infty}\sup_{t\leq T}\sup_{\vert \lambda\vert\leq k}\vert\langle X_N(t),e^{i\lambda\cdot}\rangle-\langle X_{N'}(t),e^{i\lambda\cdot}\rangle\vert=0, \quad \forall k>0 \quad \mathbb{P}-a.s.,
$$
expression (\ref{sys_3.8}) is equal to
\be
\label{sys_3.10}
\frac{2}{k}\lim_{N,N'\to\infty}\mathbb{E}\left[\int_0^T\int_{\R^d}\vert\lambda\vert^2\vert \widetilde{h_N}(\lambda,t)\vert^2+\vert\lambda\vert^2\vert\widetilde{h_{N'}}(\lambda,t)\vert^2 d\lambda dt\right].
\ee
Now
\begin{eqnarray}
\label{sys_3.11}
&&\mathbb{E}\left[\int_0^T\int_{\R^d}\vert\lambda\vert^2\vert \widetilde{h_N}(\lambda,t)\vert^2+\vert\lambda\vert^2\vert\widetilde{h_{N'}}(\lambda,t)\vert^2 d\lambda dt\right]\nonumber\\
&=&\mathbb{E}\left[\int_0^T\vert\vert \nabla h_N(\cdot,t)\vert\vert^2 dt+\int_0^T\vert\vert \nabla h_{N'}(\cdot,t)\vert\vert^2 dt\right];
\end{eqnarray}
by (\ref{eq_45}), with $A_N(t)$ as defined in (\ref{A_N}) and $S_N(t)$ defined in Proposition \ref{prop_4.4}, we
obtain that \be \mathbb{E}\left[\vert\vert \nabla
h_{N}(\cdot,t)\vert\vert^2\right]\leq
\frac{\mathbb{E}\left[S_N(T)\right]}{\sigma_N^2}+\frac{cT}{\sigma_N^2}<\infty
\ee uniformly in $N$ with $c$ positive constant. As a consequence
(\ref{sys_3.11}) is finite.

It follows that, for $k$ sufficiently large, (\ref{sys_3.10}) can
be made smaller than any given $\epsilon>0$ and there exists a
positive function $\rho(x,t)\in L^2(\R^d,\R_+)$ satisfying
equation (\ref{density}).

Since by (\ref{density}) and
$\lim_{N\to\infty}W_N(\cdot)=\delta_0$ (in the sense of
distributions)
$$
\lim_{N\to\infty}\int_{\R^d}f(x)X_N(t)(dx)=\int_{\R^d}f(x)\rho(x,t)dx\quad
f\in C_b^0(\R^d\times [0,T])\quad \mathbb{P}-a.s.,
$$
we have by (\ref{eq_3.6})
$$
\int_{\R^d}f(x,t)X(t)(dx)=\int_{\R^d}f(x,t)\rho(x,t)dx \quad f\in C_b^0(\R^d\times [0,T]), \mathbb{P}-a.s.
$$
Therefore the measure $X(t)$ is absolutely continuous with respect to the Lebesgue measure with density $\rho(x,t)$.

\hspace{14 cm} $\Box$

As a consequence of Proposition (\ref{prop_3.2})
\be
\label{eq_3.5}
\lim_{N\to\infty}\langle X_N(t),f(\cdot)\rangle=\langle X(t),f(\cdot)\rangle=\int_{\R^d}f(x)\rho(x,t)dx \quad f\in C_b^0(\R^d), t\in[0,T]
\ee

As next point we need the description of the dynamics governing the time evolution of the possible limit process
$\{X(t), t\in [0,T]\}$.

\section*{A formal derivation of the continuum models}
In this section, following \cite{mor}, we characterize the limit
behavior, as $N\to\infty$, of the process $X_N$ both in the case
$\lim_{N\to\infty}\sigma_N=0$ and
$\lim_{N\to\infty}\sigma_N=\sigma_{\infty}>0$.

By taking into account expression (\ref{eq_1.8}) and by using Ito's formula we get the following weak form of the time evolution of $X_N(t)$:
\begin{eqnarray}
\label{eq_2.2}
\langle X_N(t),f(\cdot,t)\rangle &=& \langle X_N(0),f(\cdot,0)\rangle\nonumber\\
&+&\int_0^t\langle X_N(s),(\nabla G_a*X_N(s))\cdot \nabla f(\cdot,s)\rangle ds\nonumber\\
&-&\int_0^t\langle X_N(s),\nabla g_N(\cdot,s)\cdot \nabla f(\cdot,s)\rangle ds\nonumber\\
&-&\int_0^t\langle X_N(s),\nabla U(\cdot)\cdot\nabla f(\cdot,s)\rangle ds\nonumber\\
&+&\int_0^t\langle X_N(s),\frac{1}{2}\sigma_N^2\Delta f(\cdot,s)+\frac{\partial}{\partial s}f(\cdot,s)\rangle ds\nonumber\\
&+&\frac{\sigma_N}{N}\int_0^t\sum_{k=1}^N\nabla f(X_N^k(s),s)dW_k(s), \quad f\in C_b^{2,1}(\R^d\times[0, T]).\nonumber\\
\end{eqnarray}
Last term in (\ref{eq_2.2})
$$
M_N(f,t):=\frac{\sigma_N}{N}\int_0^t\sum_{k=1}^N\nabla f(X_N^k(s),s)dW_k(s)
$$
is a martingale with respect to the natural filtration of the process $\{X_N(t), t\in[0,T]\}$ and the quadratic variation
\be
\label{quadratic}
\lim_{N\to\infty}\mathbb{E}\left[\sup_{t\leq T}\vert M_N(f,t)\vert\right]^2=0
\ee
(see \cite{mor}, \cite{morale:JMB1}). This implies, in both cases, convergence to zero in probability, that is the substantial reason of the deterministic limiting behavior of the process, as $N\to\infty$, since in this limit the evolution equation of the process will not contain the Brownian noise anymore (see \cite{morale:JMB1}).

In order to derive a formal limit for the process $X_N$ also when
$\lim_{N\to\infty}\sigma_N=0$, let us assume that $X(t)$ admits
density with respect to the Lebesgue measure also in this case. As
a formal consequence of this assumption and (\ref{eq_3.5}), we get
\begin{eqnarray*}
\lim_{N\to\infty}g_N(x,t)&=&\lim_{N\to\infty}(V_N*X_N(t))(x)=\rho(x,t),\\
\lim_{N\to\infty}\nabla g_N(x,t)&=&\nabla \rho(x,t),\\
\lim_{N\to\infty}(\nabla G_a*X_N(t))(x)&=&(\nabla G_a*X(t))(x)\\
&=&\int\nabla G_a(x-y)\rho(y,t)dy, \quad x\in\R^d, t\in [0,T].
\end{eqnarray*}

Hence by applying the above limits, from (\ref{eq_2.2}) and the hypothesis $\lim_{N\to\infty}\sigma_N=\sigma_{\infty}\geq 0$, we get the following equation
\begin{eqnarray}
\label{eq_3.8}
\int_{\R^d}f(x,t)\rho(x,t)dx \!\!\!&=&\!\!\!\!\! \int_{R^d}f(x,0)\rho(x,0)dx\nonumber\\
\!\!\!&+&\!\!\!\!\!\int_0^t\!\! ds\!\!\int_{\R^d}[(\nabla G_a*\rho(\cdot,s))(x)-\nabla U(x)-\nabla\rho(x,s)]\nonumber\\
&&\quad\quad\cdot\nabla f(x,s)\rho(x,s)dx\nonumber\\
&+&\int_0^tds\int_{\R^d}\left[\frac{\partial}{\partial s}f(x,s)\rho(x,s)+\frac{\sigma^2_{\infty}}{2}\Delta f(x,s)\rho(x,s)\right]dx.\nonumber\\
\end{eqnarray}
Equation (\ref{eq_3.8}) is the weak version of the following equation for the spatial density $\rho$:
\begin{eqnarray}
\label{integro}
\frac{\partial}{\partial t}\rho(x,t)&=&\frac{\sigma_{\infty}^2}{2}\Delta \rho(x,t)+\nabla\cdot (\rho(x,t)\nabla\rho(x,t))+\nabla\cdot(\rho(x,t)\nabla U(x))\nonumber\\
&-&\nabla\cdot[\rho(x,t)(\nabla G_a*\rho(\cdot,t))(x)], \quad x\in\R^d, t\in [0,T],\nonumber\\
\rho(x,0)&=&\rho_0(x), \quad x\in\R^d.
\end{eqnarray}

Obviously if $\sigma_{\infty}=0$ the diffusive term in (\ref{eq_3.8}) and (\ref{integro}) vanishes, while if $\sigma_{\infty}>0$ the dynamics of the density is smoothed by the diffusive term. This is due to the memory of the fluctuations existing when the number of particles $N$ is finite.

\section*{Main results}
In the present section we present the main results of this chapter,
namely two theorems on the convergence of the interacting
particle system (\ref{system_complete}) to the integro-differential
equation (\ref{integro}), both for $\sigma_{\infty}=0$ and $\sigma_{\infty}>0$.

We begin with the case $\sigma_{\infty}=0$ (non-viscous case) following the approach proposed in \cite{mor} and then we move to the case $\sigma_{\infty}> 0$ (viscous case).

\subsection*{Non-viscous case}
We are not aware of general results concerning the existence of
sufficiently regular solutions $\rho$ for equation
(\ref{integro}); therefore we need the following assumption:
\begin{assumption}
\label{ass_1} For some $T\in [0,\infty)$ system (\ref{integro})
with $\sigma_{\infty}= 0$ admits a unique, nonnegative solution
$\rho\in C_b^{[(d+2)/d]+1,1}(\R^d\times [0,T])$. 
\end{assumption}
About the uniqueness of the solution of equation (\ref{integro})
without confining potential we address to \cite{morale:MAVKDM}.

Let $\sigma_{\infty}=0$ and
suppose that
\be
\label{cond_3.2}
\widetilde{\nabla W_1}\in L^{\infty}(\R^d),
\ee

\be
\label{cond_3.3}
\vert W_1(x)\vert\leq \frac{c}{1+\vert x\vert^{d+2}},  \quad {\rm if}\quad \vert x\vert\geq 1.
\ee

Consider the following assumption for the aggregation kernel $G_a(x)$ and the confining potential $U(x)$:
\be
\label{cond_3.9}
\nabla G_a, \nabla U \in C_b^{[(d+2)/d]+1}(\R^d)\cap L^1(\R^d)
\ee

As far as $\beta$ is concerned, we need to assume that one of the following conditions is satisfied:
\be
0<\beta<\frac{d}{d+2}.
\ee
\begin{center}
or
\end{center}
\be
\label{cond_3.10}
\frac{d}{d+2}\leq \beta<1 \quad {\rm and}\quad \lim_{N\to +\infty}\sigma_NN^{\beta(d+2)/d-1}=0.
\ee

Under previous hypotheses, we can prove the following theorem:
\begin{teo}
\label{risultato_main}
Assume (\ref{cond_3.2})-(\ref{cond_3.10}) and Assumption \ref{ass_1}. If
$$
\lim_{N\to\infty}\mathbb{E}\left[\vert\vert h_N(\cdot,0)-\rho_0(\cdot)\vert\vert_2^2\right]=0,
$$
then
\be
\label{risultato_3.11}
\lim_{N\to\infty}\mathbb{E}\left[\sup_{t\leq T}\vert\vert h_N(\cdot,t)-\rho(\cdot,t)\vert\vert_2^2\right]=0,
\ee
where $\rho$ is the unique solution of (\ref{integro}) with $\sigma_{\infty}=0$.
\end{teo}

\begin{cor}
\label{cor_2.2}
Equation (\ref{risultato_3.11}) implies
$$
\lim_{N\to\infty}\langle X_N(t),f\rangle=\langle X(t),f\rangle=\int f(x)\rho(x,t)dx
$$
uniformly in $t\in [0,T]$, for any $f\in C_b^1(\R^d)\cap L^2(\R^d)$.
\end{cor}

\begin{flushleft}
{\it Proof.}
\end{flushleft}
$$
\vert \langle X_N(t)-\rho(\cdot,t),f\rangle\vert\leq \vert \langle h_N(\cdot,t)-\rho(\cdot,t),f\rangle\vert+\langle X_N(t),\vert f-f*W_N\vert\rangle;
$$
from (\ref{risultato_3.11}) and Lemma \ref{lemma_3.1} we obtains our thesis.

\hspace{14 cm} $\Box$

Previous corollary state that the empirical measure $X_N(t)$ converges weakly to the density $\rho(\cdot,t)$.

\subsubsection*{Proof of Theorem \ref{risultato_main}}
We prove this result by following the same approach used in \cite{mor} for proving it in the case of system (\ref{system_complete}) without confining potential ($\gamma_1=0$).

We obtain the convergence of $h_N(\cdot,t)$ to its limit by performing the following steps:
\begin{itemize}
\item[1.] we guess the dynamics (\ref{integro}) for the limit $\rho$;
\item[2.] we try to control $\vert\vert h_N(\cdot,t)-\rho(\cdot,t)\vert\vert_2^2$ in term of its initial value $\vert\vert h_N(\cdot,0)-\rho_0(\cdot)\vert\vert_2^2$ by writing down Ito's formula for that process and by estimating the different contributions.
\end{itemize}
We have
\be
\label{eq_3.10}
\vert\vert h_N(\cdot,t)-\rho(\cdot,t)\vert\vert_2^2=\vert\vert h_N(\cdot,t)\vert\vert_2^2+\vert\vert \rho(\cdot,t)\vert\vert_2^2-2\langle \rho(\cdot,t),h_N(\cdot,t)\rangle.
\ee
From (\ref{eq_3.8}) and integration by parts, one gets:
\begin{eqnarray}
\vert\vert \rho(\cdot,t)\vert\vert_2^2&=&\vert\vert \rho(\cdot,0)\vert\vert_2^2\nonumber\\
&+&\int_0^tds\int_{\R^d}[(\nabla G_a*\rho(\cdot,s))(x)-\nabla\rho(x,s)-\nabla U(x)]\cdot\nabla\rho(x,s)\rho(x,s)dx\nonumber\\
&+&\int_0^tds\int_{\R^d}\frac{\partial}{\partial s}(\rho(x,s))\rho(x,s)dx\nonumber\\
&=&\vert\vert \rho(\cdot,0)\vert\vert_2^2\nonumber\\
&-&\int_0^t ds\int_{\R^d}\nabla\cdot[\rho(x,s)(\nabla G_a*\rho(\cdot,s))(x)]\rho(x,s)dx\nonumber\\
&+&\int_0^tds\int_{\R^d}\nabla\cdot[\nabla\rho(x,s)\rho(x,s)]\rho(x,s)dx\nonumber\\
&+&\int_0^tds\int_{\R^d}\nabla\cdot[\nabla U(x)\rho(x,s)]\rho(x,s)dx\nonumber\\
&+&\int_0^tds\int_{\R^d}\nabla\cdot[\nabla\rho(x,s)\rho(x,s)]\rho(x,s)+\nabla\cdot[\nabla U(x)\rho(x,s)]\rho(x,s)\nonumber\\
&-&\nabla\cdot[\rho(x,s)(\nabla G_a*\rho(\cdot,s))(x)]\rho(x,s)dx\nonumber
\end{eqnarray}
\begin{eqnarray}
&=&\vert\vert \rho(\cdot,0)\vert\vert_2^2\nonumber\\
&-&2\int_0^t ds\int_{\R^d}\nabla\cdot[\rho(x,s)(\nabla G_a*\rho(\cdot,s))(x)]\rho(x,s)dx\nonumber\\
&+&2\int_0^tds\int_{\R^d}\nabla\cdot[\nabla\rho(x,s)\rho(x,s)]\rho(x,s)dx\nonumber\\
&+&2\int_0^tds\int_{\R^d}\nabla\cdot[\nabla U(x)\rho(x,s)]\rho(x,s)dx\nonumber\\
&=&\vert\vert \rho(\cdot,0)\vert\vert_2^2+2\int_0^t\langle \rho(\cdot,s),(\rho(\cdot,s)*\nabla G_a)\cdot\nabla \rho(\cdot,s)\rangle ds\nonumber\\
&-&2\int_0^t\langle \rho(\cdot,s),\vert \nabla \rho(\cdot,s)\vert^2\rangle ds-2\int_0^t\langle \rho(\cdot,s), \nabla U(\cdot)\cdot\nabla\rho(\cdot,s)\rangle ds.
\end{eqnarray}
With the same computations used to obtain expression (\ref{eq_2.15}), for the first term of (\ref{eq_3.10}) one obtains
\begin{eqnarray}
\vert\vert h_N(\cdot,t)\vert\vert_2^2&=&\vert\vert h_N(\cdot,0)\vert\vert_2^2
-2\int_0^t\langle X_N(u),\vert \nabla g_N(\cdot,s)\vert^2\rangle du\nonumber\\
&-&2\int_0^t\langle X_N(u),\nabla g_N(\cdot,u)\cdot\nabla U\rangle du\nonumber\\
&+&2\int_0^t\langle X_N(u),\nabla g_N(\cdot,u)\cdot(\nabla G_a*X_N(u))\rangle du\nonumber\\
&-&\sigma_N^2\int_0^t\vert\vert \nabla h_N(\cdot,u)\vert\vert_2^2 du\nonumber\\
&+&\frac{2\sigma_N}{N}\int_0^t\sum_{k=1}^N\nabla g_N(X_N^k(s),s)dW^k(s)\nonumber\\
&+&c_1\sigma_N^2(t-s)N^{\beta(d+2)/d-1}.
\end{eqnarray}

From (\ref{eq_2.2}), (\ref{integro}), the symmetry of $W_N$, by Ito's formula and integration by parts, one obtains
\begin{eqnarray}
&&\langle \rho(\cdot,t),h_N(\cdot,t)\rangle=\langle X_N(t),\rho(\cdot,t)*W_N\rangle\nonumber\\
&&=\langle \rho(\cdot,0),h_N(\cdot,0)\rangle\nonumber\\
&&+\int_0^t\langle X_N(t),[(X_N(s)*\nabla G_a)-\nabla g_N(\cdot,s)-\nabla U(x)]
\cdot\nabla(\rho(\cdot,t)*W_N)\rangle ds\nonumber\\
&&+\frac{\sigma_N^2}{2}\int_0^t\langle X_N(s),\Delta\rho(\cdot,s)*W_N\rangle ds
+\int_0^t\langle X_N(s),\frac{\partial}{\partial s}\rho(\cdot,s)*W_N\rangle ds\nonumber\\
&&+\frac{\sigma_N}{N}\int_0^t\sum_{k=1}^N\nabla(\rho(\cdot,s)*W_N)(X_N^k(s))dW^k(s)\nonumber
\end{eqnarray}
\begin{eqnarray}
&&=\langle \rho(\cdot,0),h_N(\cdot,0)\rangle\nonumber\\
&&+\int_0^t\langle X_N(t),[(X_N(s)*\nabla G_a)-\nabla g_N(\cdot,s)-\nabla U(x)]
\cdot\nabla\rho(\cdot,t)*W_N\rangle ds\nonumber\\
&&+\frac{\sigma_N^2}{2}\int_0^t\langle X_N(s),\Delta\rho(\cdot,s)*W_N\rangle ds\nonumber\\
&&+\int_0^t\langle X_N(s),[\nabla\cdot (\rho(\cdot,s)\nabla\rho(\cdot,s))+\nabla\cdot(\rho(\cdot,s)\nabla U)\nonumber\\
&&-\nabla\cdot[\rho(\cdot,s)(\nabla G_a*\rho(\cdot,s))]]*W_N\rangle ds
+\frac{\sigma_N}{N}\int_0^t\sum_{k=1}^N\nabla(\rho(\cdot,s)*W_N)(X_N^k(s))dW^k(s)\nonumber
\end{eqnarray}
\begin{eqnarray}
&&=\langle \rho(\cdot,0),h_N(\cdot,0)\rangle\nonumber\\
&&+\int_0^t\langle X_N(t),[(X_N(s)*\nabla G_a)-\nabla g_N(\cdot,s)-\nabla U(x)]
\cdot\nabla\rho(\cdot,t)*W_N\rangle ds\nonumber\\
&&-\frac{\sigma_N^2}{2}\int_0^t\langle \nabla h_N(\cdot,s),\nabla\rho(\cdot,s)\rangle ds\nonumber\\
&&+\!\!\!\int_0^t\!\!\!\langle \nabla h_N(\cdot,s),[-\nabla \rho(\cdot,s)-\nabla U+(\nabla G_a *\rho(\cdot,s))]\rho(\cdot,s)\rangle ds\nonumber\\
&&+\frac{\sigma_N}{N}\int_0^t\sum_{k=1}^N\nabla(\rho(\cdot,s)*W_N)(X_N^k(s))dW^k(s).
\end{eqnarray}
It follows that for (\ref{eq_3.10}) one gets the following expression:
\begin{eqnarray}\label{eq_3.14}
&&\vert\vert h_N(\cdot,t)-\rho(\cdot,t)\vert\vert_2^2\nonumber\\
&&=\vert\vert h_N(\cdot,0)-\rho(\cdot,0)\vert\vert_2^2\nonumber\\
&&-2\int_0^t\left[\langle X_N(s),\nabla g_N(\cdot,s)\cdot(\nabla g_N(\cdot,s)-\nabla\rho(\cdot,s)*W_N)\rangle\right.\nonumber\\
&&+\left.\langle \rho(\cdot,s),\nabla\rho(\cdot,s)\cdot[\nabla\rho(\cdot,s)-\nabla h_N(\cdot,s)]\rangle\right] ds\nonumber\\
&&+2\int_0^t\left[\langle X_N(s),((X_N(s)*\nabla G_a)-\nabla U)[\nabla g_N(\cdot,s)-\nabla \rho(\cdot,s)*W_N]\rangle\right.\nonumber\\
&&+\left.\langle \rho(\cdot,s),((\rho(\cdot,s)*\nabla G_a) - \nabla U)[\nabla \rho(\cdot,s)-\nabla h_N(\cdot,s)]\rangle\right]ds\nonumber\\
&&-\sigma_N^2\int_0^t\left[\vert\vert \nabla h_N(\cdot,s)\vert\vert^2_2 -\langle \nabla h_N(\cdot,s),\nabla\rho(\cdot,s)\rangle\right]ds+c_1\sigma_N^2(t-s)N^{\beta(d+2)/d-1}\nonumber\\
&&+2\frac{\sigma_N}{N}\int_0^t\sum_{k=1}^N(\nabla g_N(\cdot,s)-\nabla \rho(\cdot,s)*W_N)(X_N^k(s))dW^k(s).
\end{eqnarray}

\subsubsection*{Estimates for the terms on the right side of (\ref{eq_3.14})}
In order to get estimates for the terms on the right hand side of (\ref{eq_3.14}) we need the following lemmas:

\begin{lemma}
\label{lemma_3.1}
Let $W_1$ be a symmetric density function which satisfies (\ref{cond_3.2}) and (\ref{cond_3.3}) and let $W_N$ be defined as in (\ref{cond_3.4}). Let $f, \nabla f \in C_b^1(\R^d)$. Then for any $x\in\R^d$
$$
\vert f(x)-(f*W_N)(x)\vert\leq c_2\chi_N^{-1}\vert\vert \nabla f\vert\vert_{\infty}.
$$
Furthermore if $f, \nabla f \in L^2(\R^d)$, then
$$
\vert\vert f-f*W_N\vert\vert_2^2\leq c_2\chi_N^{-2}\vert\vert \nabla f\vert\vert_2^2.
$$
The constant $c_2$ is independent of $f$.
\end{lemma}
Lemma \ref{lemma_3.1} is proved in \cite{oel3}.

\begin{lemma}
\label{lemma_3.2}
Let us suppose Assumption \ref{ass_1} and let be $L=[(d+2)/2]$ and
$$
v(x)\in C_b^{L+1}(\R^d,\R^d),
$$
$$
G(x)\in L^1(\R^d,\R).
$$
Then we have
\begin{eqnarray*}
&&\left\vert \left\langle X_N(s)-\rho(\cdot,s),\left(\nabla [h_N(\cdot,s)-\rho(\cdot,s)]*W_N\right)\cdot v\right\rangle\right\vert\\
&&\leq c_3\left(\vert\vert h_N(s)-\rho(\cdot,s)\vert\vert_2^2+\chi_N^{-2}+N^{\beta(1-2L/d)}\right),
\end{eqnarray*}
\begin{eqnarray*}
&&\left\vert \left\langle (X_N(s)-\rho(\cdot,s))*G,\left(\nabla [h_N(\cdot,s)-\rho(\cdot,s)]*W_N\right)\cdot v\right\rangle\right\vert\\
&&\leq c_4\left(\vert\vert h_N(s)-\rho(\cdot,s)\vert\vert_2^2+\chi_N^{-2}+N^{\beta(1-2L/d)}\right).
\end{eqnarray*}
\end{lemma}
For a proof of this lemma see \cite{mor}.

We begin by considering the second term on the right hand side of (\ref{eq_3.14}); since this term does not depend on the potential $U$, we can recall a result proved in \cite{mor}.
Indeed, as showed in \cite{mor}, by considering $v=\nabla \rho(\cdot,s)*W_N$, from Lemma \ref{lemma_3.1}, Lemma \ref{lemma_3.2} and Assumption \ref{ass_1}, we have

\begin{eqnarray}
A&:=&\left\vert\langle X_N(s),\nabla g_N(\cdot,s)\cdot(\nabla g_N(\cdot,s)-\nabla\rho(\cdot,s)*W_N)\rangle\right.\nonumber\\
&&\left.+\langle \rho(\cdot,s),\nabla\rho(\cdot,s)\cdot[\nabla\rho(\cdot,s)-\nabla h_N(\cdot,s)]\rangle\right\vert\nonumber\\
&\leq& c_5\vert\vert h_N(\cdot,s)-\rho(\cdot,s)\vert\vert_2^2+\chi_N^{-2}+N^{\beta(1-2L/d)}.
\end{eqnarray}

Now let us consider the term in (\ref{eq_3.14}) involving the aggregating kernel $G_a$ and the confining potential $U$
\begin{eqnarray}
B&=&\langle X_N(s),((X_N(s)*\nabla G_a)-\nabla U)\cdot[\nabla (h_N(\cdot,s)- \rho(\cdot,s))*W_N]\rangle\nonumber\\
&-&\langle \rho(\cdot,s),((\rho(\cdot,s)*\nabla G_a) - \nabla U)\cdot[-\nabla \rho(\cdot,s)+\nabla h_N(\cdot,s)]\rangle\nonumber\\
&=&\langle X_N(s)-\rho(\cdot,s),((X_N(s)*\nabla G_a)-\nabla U)\cdot[\nabla (h_N(\cdot,s)- \rho(\cdot,s))*W_N]\rangle\nonumber\\
&+&\langle \rho(\cdot,s),((X_N(s)*\nabla G_a)-\nabla U)\cdot[\nabla (h_N(\cdot,s)- \rho(\cdot,s))*W_N]\nonumber\\
&&-((\rho(\cdot,s)*\nabla G_a) - \nabla U)\cdot[-\nabla \rho(\cdot,s)+\nabla h_N(\cdot,s)]\rangle\nonumber\\
&=&B_1+B_2.
\end{eqnarray}
By (\ref{cond_3.9}) and Lemma \ref{lemma_3.2}, with $v=-\nabla U+X_N(s)*\nabla G_a$ we get
\be
\vert B_1\vert\leq c_3\left(\vert\vert h_N(s)-\rho(\cdot,s)\vert\vert_2^2+\chi_N^{-2}+N^{\beta(1-2L/d)}\right).
\ee

On the other hand
\begin{eqnarray}
B_2&=&\langle \rho(\cdot,s),((X_N(s)*\nabla G_a)-\nabla U)\cdot[\nabla (h_N(\cdot,s)- \rho(\cdot,s))*W_N]\nonumber\\
&&-((\rho(\cdot,s)*\nabla G_a) - \nabla U)\cdot[-\nabla \rho(\cdot,s)+\nabla h_N(\cdot,s)]\rangle\nonumber\\
&=&\langle \rho(\cdot,s),[(X_N(s)-\rho(\cdot,s))*\nabla G_a]\cdot [\nabla (h_N(\cdot,s)- \rho(\cdot,s))*W_N]\rangle\nonumber\\
&&+\langle \rho(\cdot,s),((\rho(\cdot,s)*\nabla G_a)-\nabla U)\nonumber\\
&&\nabla[(h_N(\cdot,s)- \rho(\cdot,s))*W_N-(h_N(\cdot,s)- \rho(\cdot,s))]\rangle\nonumber\\
&=&\langle (X_N(s)-\rho(\cdot,s))*\nabla G_a,\rho(\cdot,s) [\nabla (h_N(\cdot,s)- \rho(\cdot,s))*W_N]\rangle\nonumber\\
&&-\langle \nabla[\rho(\cdot,s)((\rho(\cdot,s)*\nabla G_a)-\nabla U)],\nonumber\\
&&(h_N(\cdot,s)- \rho(\cdot,s))*W_N-(h_N(\cdot,s)- \rho(\cdot,s))\rangle\nonumber\\
&=&B_2^1+B_2^2.
\end{eqnarray}
By Assumption \ref{ass_1}, condition (\ref{cond_3.9}) and Lemma \ref{lemma_3.2}, with $G=\nabla G_a$ and $v=\rho-\nabla U$ we get
$$
\vert B_2^1\vert\leq c_4\left(\vert\vert h_N(\cdot,s)-\rho(\cdot,s)\vert\vert_2^2+\chi_N^{-2}+N^{\beta(1-2L/d)}\right)
$$
and by taking into account also Lemma \ref{lemma_3.1}
$$
\vert B_2^2\vert\leq c_2\left(\vert\vert h_N(\cdot,s)-\rho(\cdot,s)\vert\vert_2^2+\chi_N^{-2}\right).
$$
So for the second term $B$ of (\ref{eq_3.14}) we get the following estimate:
$$
\vert B\vert\leq c_6\left(\vert\vert h_N(\cdot,s)-\rho(\cdot,s)\vert\vert_2^2+\chi_N^{-2}+N^{\beta(1-2L/d)}\right).
$$
For the third integrand on the right side of (\ref{eq_3.14}) we have
\begin{eqnarray}
C&:=&-\sigma_N^2\left[\vert\vert \nabla h_N(\cdot,s)\vert\vert_2^2-\langle \nabla h_N(\cdot,s),\nabla\rho(\cdot,s)\rangle\right]\nonumber\\
&\leq&-\sigma_N^2\left[\frac{\vert\vert \nabla h_N(\cdot,s)\vert\vert_2^2}{2}-\langle \nabla h_N(\cdot,s),\nabla\rho(\cdot,s)\rangle\right]\nonumber\\
&=&-\frac{\sigma_N^2}{2}\vert\vert \nabla h_N(\cdot,s)-\nabla\rho(\cdot,s)\vert\vert_2^2+\frac{\sigma_N^2}{2}\vert\vert \nabla\rho(\cdot,s)\vert\vert_2^2.
\end{eqnarray}

Let us consider the submartingale term in (\ref{eq_3.14}):
$$
M_N(t):=\frac{\sigma_N}{N}\left\vert\int_0^t\sum_{k=1}^N(\nabla g_N(\cdot,s)-\nabla \rho(\cdot,s)*W_N)(X_N^k(s))dW^k(s)\right\vert.
$$
As showed in \cite{mor}, by Doob's inequality,
\be
\label{eq_3.20}
\mathbb{E}\left[\sup_{t\leq T'}M_N(t)\right]^2\leq \frac{c_7\sigma_N^2\chi_N^d}{N}\mathbb{E}\left[\int_0^{T'}\vert\vert \nabla(h_N(\cdot,s)-\rho(\cdot,s))\vert\vert_2^2ds\right], \quad T'\leq T.
\ee
By collecting all contributions, (\ref{eq_3.14}) becomes
\begin{eqnarray}
\label{eq_3.21}
\vert\vert h_N(\cdot,t)-\rho(\cdot,t)\vert\vert_2^2&\leq& \vert\vert h_N(\cdot,0)-\rho(\cdot,0)\vert\vert_2^2
+c_5\int_0^t\vert\vert h_N(\cdot,s)-\rho(\cdot,s)\vert\vert_2^2ds\nonumber\\
&+&c_5t[\chi_N^{-2}+N^{\beta(1-2L/d)}+\sigma_N^2N^{\frac{\beta(d+2)}{d}-1}+\sigma_N^2]\nonumber\\
&-&\frac{\sigma_N^2}{2}\int_0^t\vert\vert \nabla(h_N(\cdot,s)-\rho(\cdot,s))\vert\vert_2^2ds\nonumber\\
&+&\frac{2\sigma_N}{N}\left\vert\int_0^t\sum_{k=1}^N(\nabla g_N(\cdot,s)-\nabla \rho(\cdot,s)*W_N)(X_N^k(s))dW^k(s)\right\vert,\nonumber\\
\end{eqnarray}
for $0\leq t\leq T$.

From (\ref{eq_3.20}) and (\ref{eq_3.21})
\begin{eqnarray}
&&\mathbb{E}\left[\sup_{t\leq T'}\vert\vert h_N(\cdot,t)-\rho(\cdot,t)\vert\vert_2^2+\frac{\sigma_N^2}{2}\left(1-\frac{2c_7\chi_N^d}{N}\right)\int_0^{T'}\vert\vert \nabla(h_N(\cdot,t)-\rho(\cdot,t))\vert\vert_2^2dt\right]\nonumber\\
&\leq&\mathbb{E}\left[\vert\vert h_N(\cdot,0)-\rho(\cdot,0)\vert\vert_2^2\right]+c_5\int_0^{T'}\vert\vert h_N(\cdot,t)-\rho(\cdot,t)\vert\vert_2^2dt\nonumber\\
&+&c_5T'[\chi_N^{-2}+N^{\beta(1-2L/d)}+\sigma_N^2N^{\frac{\beta(d+2)}{d}-1}+\sigma_N^2].
\end{eqnarray}
For $N$ sufficiently large and by applying Gronwall's inequality we obtain
\begin{eqnarray}
&&\mathbb{E}\left[\sup_{t\leq T}\vert\vert h_N(\cdot,t)-\rho(\cdot,t)\vert\vert_2^2\right]\nonumber\\
&\leq&\left[\mathbb{E}\left[\vert\vert h_N(\cdot,0)-\rho(\cdot,0)\vert\vert_2^2\right]+T(\chi_N^{-2}+N^{\beta(1-2L/d)}+\sigma_N^2N^{\frac{\beta(d+2)}{d}-1}+\sigma_N^2)\right]e^{c_5T'}.\nonumber\\
\end{eqnarray}
As $N\rightarrow\infty$, by (\ref{cond_3.10}) and since
$\lim_{N\to\infty}\sigma_N^2N^{\frac{\beta(d+2)}{d}-1}=0$, we
obtain our thesis.

\hspace{14 cm} $\Box$

\subsection*{Viscous case}

Now we move to the case $\sigma_{\infty}>0$.

Due to technical difficulties (the presence of the non vanishing
term $\sigma_{\infty}>0$ ), in this case we can not carry out the
same proof as for Theorem \ref{risultato_main}. Therefore, by
following \cite{oel}, we try to control directly
$\mathbb{E}\left[\langle X(t),f\rangle - \langle
\rho(t),f\rangle\right]$, obtaining a result analogous to
Corollary \ref{cor_2.2}.

About the regularity and uniqueness of the solution of equation
(\ref{integro}) we make the following assumption
\begin{assumption}
\label{ass_2} System (\ref{integro})
with $\sigma_{\infty}> 0$ admits a unique, nonnegative solution
$\rho\in C^{2,1}(\R^d\times [0,T])$. 
\end{assumption}
The requirements of Assumption \ref{ass_2} are weaker than those
of Assumption \ref{ass_1} ($\rho\in C^{2,1}(\R^d\times [0,T])$
instead of $\rho\in C_b^{[(d+2)/d]+1,1}(\R^d\times [0,T])$), but
we need a further restriction on the function $W_1$ defined by
(\ref{eq_10}): $W_1$ must have compact support.

\begin{teo}
\label{teo_2.3}
If
\begin{itemize}
\item[i)]\be \label{condition_3.5}
\lim_{N\to\infty}\mathcal{L}(X_N(0))=\delta_{\mu_0}\quad {\rm
in}\quad \mathcal{{M}_{P}}(\mathcal{{M}_{P}}(\R^d)), \ee where
$\mu_0$ has density $\rho(x,0)$ with respect the Lebesgue measure,
\item[ii)]the parameter $\beta$ satisfies condition (\ref{sigma_2}),
\item[iii)]$W_1$ defined in (\ref{eq_10}) has compact support,
\end{itemize}
then \be \lim_{N\to\infty}\langle X_N(t),f(\cdot,t)\rangle=\langle
X(t),f(\cdot,t)\rangle=\int_{\R^d}f(x,t)\rho(x,t)dx \ee for any
$f\in C_b^{2,1}(\R^d,\R_+)$, where $\rho$
is the unique solution of (\ref{integro}) with
$\sigma_{\infty}>0$.
\end{teo}

\begin{flushleft}
{\it    Proof.}
\end{flushleft}
We have to show that
$$
\sup_{0\leq t\leq T}\vert\vert \langle
X(t),f(\cdot,t)\rangle - \int_{\R^d}f(x,t)\rho(x,t)dx\vert\vert_1 =0.
$$
Since (\ref{eq_3.8}) is the weak form of (\ref{integro}), it is sufficient to show that for any $f\in C_b^{2,1}(\R^d,\R_+)$,
\begin{eqnarray*}
&\mathbb{E}&\left[\left\vert \langle X(t),f(\cdot,t)\rangle-\langle \mu_0,f(\cdot,0)\rangle-\int_0^t\langle \rho(\cdot,s),\frac{1}{2}\sigma_{\infty}^2\Delta f(\cdot,s)+\frac{\partial}{\partial s}f(\cdot,s)\right.\right.\\
&&+[(\nabla G_a * \rho(\cdot,s))(\cdot)-\nabla U(\cdot)-\nabla\rho(\cdot,s)]\cdot\nabla f(\cdot,s)\rangle ds\left.\left.\right\vert\right]=0.
\end{eqnarray*}

For fixed $f\in C_b^{2,1}(\R^d,\R_+)$
\begin{eqnarray}
&\mathbb{E}&\left[\left\vert \langle X(t),f(\cdot,t)\rangle-\langle \mu_0,f(\cdot,0)\rangle-\int_0^t\langle \rho(\cdot,s),\frac{1}{2}\sigma_{\infty}^2\Delta f(\cdot,s)+\frac{\partial}{\partial s}f(\cdot,s)\right.\right.\nonumber\\
&&+[(\nabla G_a * \rho(\cdot,s))(\cdot)-\nabla U(\cdot)-\nabla\rho(\cdot,s)]\cdot\nabla f(\cdot,s)\rangle ds\left.\left.\right\vert\right]\nonumber
\end{eqnarray}
\begin{eqnarray}
&\leq&\mathbb{E}\left[\vert \langle X(t),f(\cdot,t)\rangle-\langle X_N(t),f(\cdot,t)\rangle\vert\right]\nonumber\\
&+&\mathbb{E}\left[\vert \langle \mu_0,f(\cdot,0)\rangle-\langle X_N(0),f(\cdot,0)\rangle\vert\right]\nonumber\\
&+&\frac{\sigma^2_{\infty}}{2}\mathbb{E}\left[\int_0^t\vert \langle -\rho(\cdot,s),\Delta f(\cdot,s)\rangle + \langle X_N(s),\Delta f(\cdot,s)\rangle\vert ds\right]\nonumber\\
&+&\mathbb{E}\left[\int_0^t\vert -\langle \rho(\cdot,s),\frac{\partial}{\partial s} f(\cdot,s)\rangle + \langle X_N(s),\frac{\partial}{\partial s} f(\cdot,s)\rangle\vert ds\right]\nonumber\\
&+&\mathbb{E}\left[\int_0^t\vert \langle \rho(\cdot,s),\nabla\rho(\cdot,s)\cdot\nabla f(\cdot,s)\rangle - \langle h_N(\cdot,s),\nabla h_N(\cdot,s)\cdot\nabla f(\cdot,s)\rangle\vert ds\right]\nonumber\\
&+&\mathbb{E}\left[\left\vert\int_0^t \langle h_N(\cdot,s),\nabla h_N(\cdot,s)\cdot\nabla f(\cdot,s)\rangle - \langle X_N(s),\nabla g_N\cdot\nabla f(\cdot,s)\rangle ds\right\vert\right]\nonumber\\
&+&\mathbb{E}\left[\int_0^t\left\vert -\langle \rho(\cdot,s),[(\nabla G_a*\rho(\cdot,s))(\cdot)-\nabla U(\cdot)]\cdot\nabla f(\cdot,s)\rangle \right.\right.\nonumber\\
&&+\left.\left. \langle X_N(s),[(\nabla G_a*X_N(s))(\cdot)-\nabla U(\cdot)]\cdot\nabla f(\cdot,s)\rangle\right\vert ds\right]\nonumber\\
&+&\mathbb{E}\left[\left\vert \frac{\sigma_N}{N}\int_0^t\sum_{k=1}^N\nabla f(X_N^k(s),s)dW_k(s)\right\vert\right]\nonumber\\
&+&\mathbb{E}\left[\left\vert \langle X_N(t),f(\cdot,t)\rangle - \langle X_N(0),f(\cdot,0)\rangle \right.\right.\nonumber\\
&&-\left.\left. \int_0^t\langle X_N(s),(\nabla G_a*X_N(s))\cdot \nabla f(\cdot,s)\rangle ds +\int_0^t\langle X_N(s),\nabla g_N(\cdot,s)\cdot \nabla f(\cdot,s)\rangle ds \right.\right.\nonumber\\
&&+\left.\left. \int_0^t\langle X_N(s),\nabla U(\cdot)\cdot\nabla f(\cdot,s)\rangle ds - \int_0^t\langle X_N(s),\frac{1}{2}\sigma_N^2\Delta f(\cdot,s)+\frac{\partial}{\partial s}f(\cdot,s)\rangle ds \right.\right.\nonumber\\  &&-\left.\left. \frac{\sigma_N}{N}\int_0^t\sum_{k=1}^N\nabla f(X_N^k(s),s)dW_k(s) \right\vert\right]\nonumber\\
&:=&\sum_{i=1}^9I_N^i(t).
\end{eqnarray}

Clearly, by (\ref{eq_3.5}) and hypothesis (\ref{condition_3.5}),
$\sum_{i=1}^4I_N^i(t)=0$, by (\ref{eq_2.2}) $I_N^9(t)=0$ and by
(\ref{quadratic}) $I_N^8(t)=0$. It remains to estimate the terms
$I_N^5(t)$, $I_N^6(t)$ and $I_N^7(t)$.
\begin{eqnarray*}
I_N^5(t)&=&\mathbb{E}\left[\int_0^t\vert \langle \rho(\cdot,s),\rho(\cdot,s)\Delta f(\cdot,s)\rangle - \langle h_N(\cdot,s), h_N(\cdot,s)\Delta f(\cdot,s)\rangle\vert ds\right]\\
&\leq&\vert\vert \Delta f\vert\vert_{\infty}\int_0^T\mathbb{E}\left[\int_{\R^d}\vert h_N(x,t)-\rho(x,t)\vert\vert h_N(x,t)+\rho(x,t)\vert dx\right]dt\\
&\leq&\vert\vert \Delta f\vert\vert_{\infty}\left(\mathbb{E}\left[\int_0^T\int_{\R^d}\vert h_N(x,t)-\rho(x,t)\vert^2dx dt\right]\right)^{1/2}\\
&&\cdot\left(\mathbb{E}\left[\int_0^T\int_{\R^d}\vert h_N(x,t)+\rho(x,t)\vert^2dx dt\right]\right)^{1/2};
\end{eqnarray*}
by (\ref{eq_45}) and (\ref{density}) we obtain \be
\lim_{N\to\infty}I_N^5(t)=0. \ee By the symmetry of $W_1$,
\begin{eqnarray}
\label{I6} I_N^6(t)&=&\mathbb{E}\left[\left\vert \int_0^t \langle
X_N(s),W_N*(\nabla h_N(\cdot,s)\cdot\nabla
f(\cdot,s))\right.\right.\nonumber\\
&&\left.\left.-(W_N*\nabla h_N(\cdot,s))\cdot\nabla f(\cdot,s)\rangle ds \right\vert\right]\nonumber\\
&=&\mathbb{E}\left[\left\vert \int_0^t
\left(\int_{\R^d}X_N(s)(dx)\int_{\R^d}W_N(x-y)\nabla
h_N(y,s)\right.\right.\right.\nonumber\\
&&\left.\left.\left.\cdot (\nabla f(y)-\nabla f(x))dy\right)ds \right\vert\right].\nonumber\\
\end{eqnarray}
By the definition of $W_N$ and since $W_1$ has compact support, with $c={\rm diam}({\rm supp}W_1(\cdot))$ and $\vert\vert D^2f\vert\vert_{\infty}=\sup_{i,j\leq d}\vert\vert \partial_{ij}^2\vert\vert_{\infty}$, (\ref{I6}) is less than or equal to
\begin{eqnarray*}
&&c\chi_N^{-1}\vert\vert D^2f\vert\vert_{\infty}\mathbb{E}\left[\int_0^t\langle X_N(s)*W_N,\vert \nabla h_N(\cdot,s)\vert\rangle ds\right]\\
&\leq&c\chi_N^{-1}\vert\vert D^2f\vert\vert_{\infty}\left(\mathbb{E}\left[\int_0^T\vert\vert h_N(\cdot,s)\vert\vert_2^2ds\right]\right)^{1/2}\left(\mathbb{E}\left[\int_0^T\vert\vert \nabla h_N(\cdot,s)\vert\vert_2^2ds\right]\right)^{1/2}\\
&\stackrel{(\ref{eq_45})}{\leq}&c\chi_N^{-1}\vert\vert D^2f\vert\vert_{\infty}.
\end{eqnarray*}
It follows that \be \lim_{N\to\infty}I_N^6(t)=0. \ee

\begin{eqnarray*}
I_N^7(t)&=&\mathbb{E}\left[\int_0^t\left\vert -\langle \rho(\cdot,s),[(\nabla G_a*\rho(\cdot,s))(\cdot)-\nabla U(\cdot)]\cdot\nabla f(\cdot,s)\rangle \right.\right.\nonumber\\
&&+\left.\left. \langle X_N(s),[(\nabla G_a*X_N(s))(\cdot)-\nabla U(\cdot)]\cdot\nabla f(\cdot,s)\rangle\right.\right.\nonumber\\
&&+\langle X_N(s),[(\nabla G_a*\rho(\cdot,s))(\cdot)-\nabla
U(\cdot)]\cdot\nabla
f(\cdot,s)\rangle\\
&&-\left.\left. \langle X_N(s),[(\nabla
G_a*\rho(\cdot,s))(\cdot)-\nabla U(\cdot)]\cdot\nabla
f(\cdot,s)\rangle\right\vert ds\right]\\
&\leq&\mathbb{E}\left[\int_0^t\left\vert \langle
X_N(s)-\rho(\cdot,s),[(\nabla G_a*\rho(\cdot,s))(\cdot)-\nabla
U(\cdot)]\cdot\nabla f(\cdot,s)\rangle\right.\vert\right.\\
&&+\left.\left\vert  \langle X_N(s),[(\nabla
G_a*\rho(\cdot,s))(\cdot)-(\nabla G_a*X_N(s))(\cdot)]\cdot\nabla
f(\cdot,s)\rangle \right\vert ds\right].
\end{eqnarray*}
By Assumption \ref{ass_2} and (\ref{eq_3.6})
$$
\lim_{N\to\infty}I_N^7(t)=0.
$$

As a consequence
$$
\lim_{N\to\infty}\sum_{i=1}^9I_N^i(t)=0.
$$

\hspace{14 cm} $\Box$

\section*{Conclusions}
In this paper we have studied the asymptotic behavior of system \eqref{system_complete}
for the size of the population $N$ growing to infinity,
being the time $t$ fixed, in terms of a law of large numbers for the empirical process $\{X_N(t), t\in\R_+\}$.

It is also of interest to study the limiting behavior of such a system for fixed $N$  and time growing to infinity. In \cite{ortisi} the authors investigate conditions for the existence of an invariant measure for system \eqref{system_complete},
i.e. conditions about the interaction potential and the confining potential such that there exists an invariant measure for the particle positions and, as a consequence, for the empirical process $\{X_N(t), t\in\R_+\}$.

\end{document}